\documentclass[11pt]{article}

\usepackage{amsmath,amssymb,amsfonts,textcomp,amsthm,xifthen,psfrag,graphicx,color,MnSymbol}
\usepackage[T1]{fontenc}
\usepackage{hyperref}
\usepackage{placeins}

\usepackage{pgfplots}
\pgfplotsset{compat=newest}
\usepgfplotslibrary{groupplots}
\usepgfplotslibrary{external}
\definecolor{navyblue}{rgb}{0.0, 0.0, 0.5}

\date{}

\newtheorem{theorem}{Theorem}
\newtheorem{lemma}[theorem]{Lemma}

\newtheorem{prop}[theorem]{Proposition}
\newtheorem{remark}[theorem]{Remark}
\theoremstyle{definition} 

\newcommand{\ndof}{\#\mathrm{dof}}

\newcommand{\<}{\langle{}}
\renewcommand{\>}{\rangle}

\newcommand{\ip}[2]{\llangle#1\hspace*{.5mm},#2\rrangle}
\newcommand{\dual}[2]{\<#1\hspace*{.5mm},#2\>}

\newcommand{\vdual}[2]{(#1\hspace*{.5mm},#2)}

\renewcommand{\skew}[1]{\mathrm{skew}(#1)}
\newcommand{\diam}{\mathrm{diam}}

\newcommand{\wat}{\widehat}

\newcommand{\cQ}{{c_Q}}
\newcommand{\cQinv}{{c_Q^{-1}}}
\newcommand{\Cdisp}{{\mathbf{C}_\mathrm{disp}}}
\newcommand{\Cdispt}{{\mathbf{C}_\mathrm{disp}^2}}
\newcommand{\Cdispinv}{{\mathbf{C}_\mathrm{disp}^{-1}}}

\newcommand{\Grad}{\boldsymbol{\nabla}}
\newcommand{\sGrad}{\boldsymbol{\varepsilon}}
\newcommand{\kGrad}{{\boldsymbol{\nabla}'}}

\def\Div{{\rm\bf div\,}}

\def\grad{\nabla}

\def\GG{\mathbf{G}}
\def\MM{\mathbf{M}}
\def\SS{\mathbf{S}}
\def\BB{\mathbf{B}}
\def\NN{\mathbf{N}}
\def\TT{\mathbf{T}}

\def\QQ{\mathbf{Q}}

\newcommand{\bL}{\ensuremath{\mathbf{L}}}
\newcommand{\LL}{\ensuremath{\mathbb{L}}}
\def\tq{\wat{\boldsymbol{q}}}

\def\bbeta{\boldsymbol{\beta}}
\def\bkappa{\boldsymbol{\kappa}}

\def\tMM{\wat{\mathrm{M}}}

\def\tNNn{\wat{\mathrm{N}}_n}

\def\tbu{\wat{\boldsymbol{u}}}

\def\tw{\wat{w}}

\newcommand{\bg}{\boldsymbol{g}}

\newcommand{\bp}{\boldsymbol{p}}
\newcommand{\bv}{\boldsymbol{v}}
\newcommand{\bu}{\boldsymbol{u}}

\newcommand{\uu}{\mathfrak{u}}
\newcommand{\tuu}{\wat{\mathfrak{u}}}

\newcommand{\deltavv}{\delta\!\vv}
\newcommand{\deltauu}{\delta\!\uu}
\newcommand{\vv}{\mathfrak{v}}
\newcommand{\ww}{\mathfrak{w}}
\newcommand{\UBC}[1]{{U_{#1}}}
\newcommand{\UU}{\ensuremath{\mathfrak{U}}}
\newcommand{\VV}{\ensuremath{\mathfrak{V}}}
\newcommand{\VVBC}[1]{\ensuremath{\mathfrak{V}_{#1}}}

\newcommand{\deltatuu}{\delta\!\wat{\mathfrak{u}}}
\newcommand{\deltaw}{\delta\!w}
\newcommand{\deltaz}{\delta\!z}

\newcommand{\deltabu}{{\boldsymbol{\delta}\!\bu}}
\newcommand{\deltabv}{{\boldsymbol{\delta}\!\bv}}

\newcommand{\deltaMM}{\boldsymbol{\delta}\!\mathbf{M}}

\newcommand{\deltaNN}{\boldsymbol{\delta}\!\mathbf{N}}

\newcommand{\traceDD}[1]{\mathrm{tr}_{#1}^{\mathrm{dDiv}}} 
\newcommand{\traceGG}[1]{\mathrm{tr}_{#1}^{\mathrm{Ggrad}}} 
\newcommand{\traceSK}[1]{\mathrm{tr}_{#1}^{\mathrm{sK}}}
\newcommand{\traceG}[1]{\mathrm{tr}_{#1}^{\mathrm{Grad}}}
\newcommand{\traceD}[1]{\mathrm{tr}_{#1}^{\mathrm{Div}}}

\newcommand{\dDiv}{{\div\Div\!}}
\newcommand{\Hdiv}[1]{{H(\div\!,#1)}}
\newcommand{\HDiv}[1]{{\HH(\Div\!,#1)}}
\newcommand{\HDivs}[1]{{\HH^s(\Div\!,#1)}}

\newcommand{\HdDiv}[1]{{\HH(\dDiv,#1)}}
\newcommand{\HdDivz}[1]{{\HH_0(\dDiv,#1)}}

\newcommand{\HSK}{\ensuremath{\mathbf{H}^{\mathrm{sK}}}}
\newcommand{\HSKBC}[1]{\ensuremath{\mathbf{H}^{\mathrm{sK}}_{#1}}}
\newcommand{\HH}{\ensuremath{\mathbb{H}}}
\newcommand{\bH}{\ensuremath{\mathbf{H}}}

\newcommand{\err}{\mathrm{err}}

\newcommand{\trSK}[1]{{\mathrm{sK},#1}}

\def\div{{\rm div\,}}

\newcommand{\ttt}{{\mathfrak{T}}}
\def\tangent{\bt}

\newcommand{\R}{\ensuremath{\mathbb{R}}}
\newcommand{\N}{\ensuremath{\mathbb{N}}}

\newcommand{\nn}{\ensuremath{\mathbf{n}}}

\newcommand{\cC}{\ensuremath{\mathcal{C}}}

\newcommand{\cT}{\ensuremath{\mathcal{T}}}

\newcommand{\cS}{\ensuremath{\mathcal{S}}}

\newcommand{\bt}{\ensuremath{\mathbf{t}}}
\newcommand{\cP}{\ensuremath{\mathcal{P}}}
\newcommand{\OO}{\ensuremath{\mathcal{O}}}

\newcommand{\xx}{{\boldsymbol{x}}}

\title{A DPG method for shallow shells
\thanks{Supported by ANID-Chile through FONDECYT projects 1190009, 1210391,
        Ruth och Nils-Erik Stenb\"{a}cks stiftelse and Oulun rakennustekniikan s\"{a}\"{a}ti\"{o}}}
\author{
Thomas~F\"uhrer$^\dagger$
\and
Norbert Heuer\thanks{
Facultad de Matem\'aticas, Pontificia Universidad Cat\'olica de Chile,
Avenida Vicu\~na Mackenna 4860, Santiago, Chile,
email: {\tt \{tofuhrer,nheuer\}@mat.uc.cl}}
\and
Antti H. Niemi\thanks{
Civil Engineering Research Unit, Faculty of Technology, University of Oulu,
Erkki Koiso-Kanttilan katu 5, 90570 Oulu, Finland,
email: {\tt antti.niemi@oulu.fi}}}

\begin{document}
\maketitle
\begin{abstract}
We develop and analyze a discontinuous Petrov--Galerkin method with optimal test functions
(DPG method) for a shallow shell model of Koiter type. It is based on a uniformly
stable ultraweak formulation and thus converges robustly quasi-uniformly.
Numerical experiments for various cases, including the Scordelis--Lo cylindrical roof,
elliptic and hyperbolic geometries, illustrate its performance. The built-in DPG error estimator gives
rise to adaptive mesh refinements that are capable to resolve boundary and interior layers.
The membrane locking is dealt with by raising the polynomial degree only of the
tangential displacement trace variable.

\end{abstract}

\section{Introduction}

The purpose of this paper is to design a quasi-optimal Galerkin scheme for a shallow
shell model of Koiter type that is uniformly stable with respect to important parameters
like the shell thickness and its curvature. Our scheme has a built-in error estimator
comprised of local contributions. It can be used to steer adaptive mesh refinements,
to efficiently solve problems with layers and singularities, including point loads.

Our scheme is a DPG method, short for discontinuous Petrov--Galerkin method with optimal test functions.
This is a Galerkin framework proposed by Demkowicz and Gopalakrishnan
\cite{DemkowiczG_10_CDP,DemkowiczG_11_CDP,DemkowiczGN_12_CDP}.
It aims at an automatic discrete inf-sup stability and
usually combines an ultraweak variational formulation \cite{Despres_94_SUF,CessenatD_98_AUW}
and independent trace variables \cite{BottassoMS_02_DPG}
with the use of optimal test functions \cite{BarrettM_84_ASP} from product Hilbert spaces.
For an early overview and its relation with mixed formulations and the minimum residual
method, see \cite{DemkowiczG_14_ODM}.

Recently we have exploited the flexibility of the DPG framework to design quasi-optimal
Galerkin schemes for the Kirchhoff--Love plate bending model, including the case
of non-convex polygonal plates \cite{FuehrerHN_19_UFK,FuehrerH_19_FDD,FuehrerHH_21_TOB}.
In this paper we extend the ultraweak variational formulation and the DPG method developed in
\cite{FuehrerHN_19_UFK} to a shallow shell model of Koiter type which has its roots in the vast
literature on shell theory.
The shallow shell model has been formulated rather concisely in the celebrated textbook by Fl\"{u}gge
\cite{Fluegge_60_SS}. It was revived by Pitk\"{a}ranta {\it et~al.} in \cite{PitkaerantaMS_01_FMA}
where also a Naghdi type variant with non vanishing transverse shear deformations was studied.
For more information on the underlying shell theory, we refer the reader to Piila who,
also together with Pitk\"{a}ranta, provided asymptotic analyses for different
geometries~\cite{Piila_94_CMT,Piila_96_CMT,PiilaP_93_CMT,PiilaP_93_EER,PiilaP_95_EER},
see also \cite{PitkaerantaLOP_95_SDS,RammW_04_SSS}.
More recently, numerical analysis for the general Koiter model has been provided in
\cite{BechetSPM_10_SST,BenzakenEMcCT_21_NML,MerabetN_17_PML,RafetsederZ_19_NMA,Zhang_16_ADG}.

Restricting the discussion to shells which deviate only slightly from flat plates substantially
reduces the technical difficulties of the general shell theory. At the same time, the shallow shell model
allows treatment of various relevant engineering problems featuring different shell geometries
(parabolic, elliptic, hyperbolic) and preserves the main characteristics of shell mechanics
from the numerical analysis point of view. It should also be noted that any smooth shell is
shallow in the vicinity of a given point and may be divided into parts that can be analyzed
by the shallow shell theory, cf.~\cite{Malinen_01_OCS,Niemi_10_BSE,Niemi_16_BCS}

A properly designed DPG method combines the advantages of mixed and
least squares finite elements by providing direct approximations of both stress and displacement
variables as well as a built-in residual for a posteriori error control,
see~\cite{DemkowiczGN_12_CDP,CarstensenDG_14_PEC}.
Our present method features direct approximation of the most important quantities of interest
for the shell problem, i.e., normal forces and bending moments together with the
middle surface displacements. The evaluation of transverse shear forces and normal rotations still
requires post processing.
The performance and robustness of our method in approximating different shell deformation states
is demonstrated in a variety of benchmark problems involving different shell geometries.
In particular, the ability of the adaptive algorithm to capture boundary and interior layers
together with the associated stress concentrations is demonstrated. The problem of membrane locking
is dealt with by utilizing a higher degree of approximation for the
variable that represents the trace of the tangential displacements.
As expected, the problem is most severe for pure bending deformations
with vanishing membrane strains but it is present in a milder form for any deformation producing
bending energy like the edge effects, see, e.g.,
\cite{Pitkaeranta_92_PML,ChapelleB_98_FCF,PitkaerantaMS_01_FMA,Niemi_08_ASL}.

Let us comment on some mathematical aspects of our approach for the shallow shell model.
Central point for the analysis of ultraweak formulations with product test spaces, as used here,
is the study of corresponding trace operators and their images.
It is critical to find a test norm that ``characterizes'' the variational formulation,
unknown at that stage, in the sense that it ensures uniform stability of the adjoint problem and
uniform boundedness of the bilinear form of the variational formulation.
There is a certain canonical procedure based on formal adjoint operators and graph norms,
cf.~\cite{CarstensenDG_16_BSF} and \cite[Appendix A]{DemkowiczGNS_17_SDM}.
For singularly perturbed problems, however, the whole setup has to be chosen very carefully,
cf.~\cite{DemkowiczH_13_RDM} for convection-dominated diffusion,
\cite{HeuerK_17_RDM} for reaction-dominated diffusion, and
\cite{FuehrerHS_20_UFR} for the Reissner--Mindlin plate bending model.
In the case of the shallow shell model, key issue is to introduce norms that take model parameters
into account, including the shell thickness, its curvature, the size of the domain, and the type and
location of boundary conditions. In \cite{FuehrerH_21_RDM} we have seen that already
moderately sized domains may cause a locking phenomenon when test norms ignore the extension of the
domain. Let us note for that matter that we incorporate the size of the domain by the parameter $D$
and the boundary conditions through the selection of the tensor $\Cdisp$,
cf.~\eqref{Poincare} below.

In practice, optimal test functions of DPG schemes have to be approximated. Theoretical results
carry over to this case once the existence of a Fortin operator is established,
cf.~\cite{GopalakrishnanQ_14_APD}. For the terms appearing in our formulation, Fortin operators
have been studied in \cite{GopalakrishnanQ_14_APD,FuehrerH_19_FDD}. But at this point it is open
how to construct variants for our Koiter model that are uniformly bounded with respect to critical
model parameters. We do not consider this problem here.

An overview of the remainder of this paper is as follows. In the next section we
introduce the model problem of a shallow Koiter shell, present a scaled formulation
and briefly discuss boundary conditions.
Section~\ref{sec_traces} is devoted to the definition and analysis of trace operators
and norms. Principal results in this section are Propositions~\ref{prop_norm_sK}
and~\ref{prop_jumps_sK}, respectively proving the identity of two trace norms and
characterizing the continuity of test functions.
However, as noted before, main contribution is the proper selection of spaces and norms.
Our ultraweak formulation is presented in Section~\ref{sec_DPG}.
Furthermore, we claim its uniform stability (Theorem~\ref{thm_stab}),
define the DPG approximation scheme and state its uniform quasi-optimal convergence
(Theorem~\ref{thm_DPG}). Proofs of Theorems~\ref{thm_stab} and~\ref{thm_DPG} are given
in Section~\ref{sec_adj}. In Section~\ref{sec_num} we report on various numerical examples
that underline the performance of our DPG scheme.
We consider shells with different geometries and show that the built-in error estimator
generates appropriate mesh refinements. In particular, it detects interior and boundary layers.

Throughout the paper, $a\lesssim b$ means that $a\le cb$ with a generic constant $c>0$
that is independent of important parameters: the (scaled) shell thickness $d$,
the diameter of $\Omega$ (the parameter domain in $\R^2$ of the mid-surface),
the curvature tensor $\BB$, and the underlying mesh $\cT$. Similarly, we use the notation $a\gtrsim b$.

\section{Model problem} \label{sec_model}

Let $\Omega\subset\R^2$ be a bounded, simply connected Lipschitz domain
with boundary $\Gamma=\partial\Omega$.
For a vertical load $f$ and a tangential load $\bp$ we consider
\begin{align} \label{prob_t}
     \BB\vdotdot\NN -\div\Div\MM &= f,\quad
    \MM - \cC_b\bkappa = 0,\quad
    \NN - \cC_m\bbeta  =0,\quad
    -\Div\NN = \bp\quad\text{in}\quad\Omega.
\end{align}
Here,
\[
   \bbeta = \sGrad{\bu}+\BB w,\quad \bkappa = -\sGrad{\grad w},\quad
   \cC_b = \frac {d^3 E}{12}\cC,\quad
   \cC_m = d E \cC,
\]
$\bbeta$ are the membrane strains, $\bu=(u_1,u_2)$ the tangential displacements,
$w$ is the transverse deflection,
$\NN$ the membrane forces, $\bkappa$ the bending curvatures,
$\MM$ the bending moments, and $\BB=\BB^T$ the geometric curvature tensor of the shell mid-surface.
The operator $\sGrad$ is the symmetric gradient,
$\sGrad\bu:=\sGrad(\bu):=\frac 12(\Grad\bu + \Grad\bu^T)$.
We assume isotropic, linearly elastic material so that $\cC$ is a positive definite,
symmetric fourth-order tensor
incorporating the plane stress assumption depending on the Poisson ratio $\nu$,
whereas we have written the dependence on $d$ and the Young modulus $E$ explicitly.

In the following we consider the rescaled displacement $dE\bu\to\bu$ and deflection $dEw\to w$.
Then, problem \eqref{prob_t} becomes
\begin{subequations} \label{prob}
\begin{alignat}{2}
     \BB\vdotdot\NN -\div\Div\MM &= f  && \quad\text{in} \quad \Omega\label{p1},\\
    \MM + \frac{d^2}{12} \cC \sGrad{\grad w} &= 0  && \quad\text{in} \quad \Omega\label{p2},\\
    \NN - \cC(\sGrad{\bu}+\BB w)  &=0   && \quad\text{in} \quad \Omega\label{p3},\\
    -\Div\NN &= \bp && \quad\text{in} \quad \Omega\label{p4}.
\end{alignat}
\end{subequations}
It remains to specify boundary conditions.
For ease of presentation we only consider two cases, a simply supported shell,
\begin{equation} \label{BCs}
   \bu=0,\quad w=0,\quad\nn\cdot\MM\nn=0\qquad\text{on}\quad \Gamma
\end{equation}
and a clamped shell,
\begin{equation} \label{BCc}
   \bu=0,\quad w=0,\quad \partial_\nn w=0 \qquad\text{on}\quad \Gamma.
\end{equation}
Here, $\nn$ denotes the exterior unit normal vector along $\Gamma$,
and $\partial_\nn w$ is the exterior normal derivative of $w$.
However, for a piecewise smooth boundary $\Gamma$,
our theory can be easily extended to mixed cases and any combination
of physically meaningful conditions including effective shear forces and point
loads. This is due to our mathematical formulation with trace operators that are
well defined for polygonal domains (and extend to piecewise smooth boundaries),
cf. our setting for Kirchhoff--Love plates \cite{FuehrerHN_19_UFK}.
The (natural) condition is that there are sufficient restrictions on the displacement
and the deflection to ensure their uniqueness. This enters through the validity of appropriate
Poincar\'e-type inequalities. It is well known that such inequalities do not hold uniformly
with respect to the diameter of the underlying physical domain. More importantly,
when ignoring such a dependence in the design of minimum residual methods in general,
and the DPG method in particular, approximations suffer from locking already
for domains of moderate diameter, see~\cite{FuehrerH_21_RDM}.
We therefore follow the strategy from \cite{FuehrerH_21_RDM} which consists
in employing scaled norms with parameters that ensure the uniform validity of 
required Poincar\'e estimates. In the case of the shallow shell of Koiter type we assume
the existence of a symmetric, positive definite tensor $\Cdisp$ of order two and a
positive number $D$ such that
\begin{equation} \label{Poincare}
   D^{-2}\|\Cdisp \bu\|^2
   \lesssim \|\sGrad{\bu}-\BB w\|^2 +  d^2 \|\sGrad{\grad w}\|^2,\quad
   \|w\|^2 \lesssim  D^4 \|\sGrad{\grad w}\|^2
\end{equation}
hold with generic constants that are independent of $d$, $\diam(\Omega)$, $\BB$,
and the involved functions. In this paper, $\Cdisp$ and $D$ are constant.
But estimates can be easily generalized to $L_\infty$ tensor/scalar functions $\Cdisp$, $D$
with corresponding uniform properties.
In the case of boundary conditions \eqref{BCs} and \eqref{BCc}, it is easy to see that
\[
   D:=\diam(\Omega),\quad
   \Cdisp:=\mathrm{diag}(c,c)\in\R^{2\times 2}
   \quad\text{with}\quad
   c := \min\Bigl\{1,\frac {d}{D^2\|\BB\|_\infty}\Bigr\},\
   \min\{1,\frac 10\}:=1
\]
are possible choices for \eqref{Poincare} to hold.
However, for a stretched domain, $D$ can be chosen as the minimum width
of $\Omega$. In general, the best choice for $D$ and $\Cdisp$ depends on the
boundary conditions, the shape of $\Omega$ and the geometric curvature tensor $\BB$,
and in some cases $\Cdisp$ reflects an anisotropy in the stability of the tangential displacements.
We refer to Section~\ref{sec_num} for specific examples where better (stronger)
selections of $\Cdisp$ are possible.
For our analysis we need the Poincar\'e estimate \eqref{Poincare} for test functions, cf.~\eqref{ass}.

\section{Trace spaces and norms} \label{sec_traces}

In the following we present and study several trace operators, product variants with support
on $\cS$ and domain variants with support on $\Gamma$. We first recall canonical
trace operators related with differential operators of order one. Then, in \S\ref{sec_KL},
we recall notation and results for operators of second order. Finally, in \S\ref{sec_sK},
we introduce a combination of trace operators that characterizes traces induced by an
ultraweak formulation of problem~\eqref{prob}.

Before doing so let us fix some notation.
Let $\cT$ be a mesh consisting of general non-intersecting Lipschitz elements $\{T\}$
so that $\bar\Omega=\cup\{\bar T;\; T\in\cT\}$. We formally denote the mesh skeleton
by $\cS=\{\partial T;\;T\in\cT\}$. For a subdomain $\omega\subset\Omega$, generally
$\omega\in\cT$ or $\omega=\Omega$, we use the standard $L_2$ spaces
$L_2(\omega)$, $\bL_2(\omega)=(L_2(\Omega))^2$ and $\LL_2(\omega)=(L_2(\Omega))^{2\times 2}$
of scalar, vector and tensor fields of order $2$, respectively, with generic norm $\|\cdot\|_\omega$
and duality $\vdual{\cdot}{\cdot}_\omega$. We drop the index $\omega$ when $\omega=\Omega$.
Furthermore, $\bH^1(\omega)=(H^1(\omega))^2$ and $H^2(\omega)$ are the standard Sobolev spaces,
$\bH^1_0(\Omega)$ is the space of $\bH^1(\Omega)$-functions with zero trace on $\Gamma$,
$H^2_0(\Omega)$ is the space of $H^2(\Omega)$-functions with zero trace and normal derivative
on $\Gamma$, $\LL_2^s(\omega),\LL_2^k(\omega)\subset\LL_2(\omega)$
are the subspaces of symmetric and skew-symmetric $L_2$-tensors, respectively, and
\begin{align*}
   &\HDiv{\omega} := \{\TT\in\LL_2(\omega);\; \Div\TT\in \bL_2(\omega)\},\qquad
   \HDivs{\omega} := \HDiv{\omega}\cap \LL_2^s(\omega),\\
   &\HdDiv{\omega} := \{\SS\in\LL_2^s(\omega);\; \div\Div\SS\in L_2(\omega)\}.
\end{align*}
We use the corresponding product spaces, denoted in the same way but replacing $\omega$ with $\cT$,
e.g., $H^1(\cT)=\Pi_{T\in\cT} H^1(T)$. The $L_2(\cT)$-dualities are denoted by $\vdual{\cdot}{\cdot}_\cT$.

\subsection{Canonical trace operators}

We define operators
\begin{align*} 
   \traceG{}:\;
   \begin{cases}
      \bH^1(\cT) &\to\ \HDivs{\cT}',\\
      \quad \bv    &\mapsto\ \dual{\traceG{}(\bv)}{\TT}_\cS
                    := \vdual{\bv}{\Div\TT}_\cT + \vdual{\sGrad{\bv}}{\TT}_\cT
   \end{cases}
\end{align*}
and
\begin{align*} 
   \traceD{}:\;
   \begin{cases}
      \HDiv{\cT} &\to \bH^1(\cT)' \\
      \quad \TT    &\mapsto\ \dual{\traceD{}(\TT)}{\bv}_\cS
                    := \vdual{\TT}{\Grad{\bv}}_\cT + \vdual{\Div\TT}{\bv}_\cT.
   \end{cases}
\end{align*}
Below, we will systematically use that $\vdual{\Grad\bv}{\TT}_\omega=\vdual{\sGrad\bv}{\TT}_\omega$
for $\bv\in\bH^1(\omega)$, $\TT\in\LL_2^s(\omega)$, $\omega\subset\Omega$. That is,
\[
   \dual{\traceG{}(\bv)}{\TT}_\cS = \dual{\traceD{}(\TT)}{\bv}_\cS
   \quad (\bv\in\bH^1(\cT),\ \TT\in\HDivs{\cT}).
\]
Analogously we define the canonical trace operators for $\Omega$ (instead of $\cT$)
with support on $\Gamma$,
$\traceG{\Omega}:\;\bH^1(\Omega)\to\HDivs{\Omega}'$ and
$\traceD{\Omega}:\;\HDiv{\Omega}\to\bH^1(\Omega)'$.

\begin{prop} \label{prop_can}
(i) If $\bv\in \bH^1(\cT)$, then
\[
   \bv\in \bH^1_0(\Omega) \quad\Leftrightarrow\quad
   \dual{\traceD{}(\TT)}{\bv}_\cS = 0
   \quad\forall \TT\in \HDiv{\Omega}.
\]
(ii) If $\TT\in \HDivs{\cT}$, then
\[
   \TT\in \HDivs{\Omega} \quad\Leftrightarrow\quad
   \dual{\traceG{}(\bv)}{\TT}_\cS = 0
   \quad\forall \bv\in \bH^1_0(\Omega).
\]
\end{prop}

\begin{proof}
For the scalar/vector cases of $H^1(\cT)$ and $\Hdiv{\cT}$, the statements
are proved in \cite[Theorem~2.3, Remark~2.5]{CarstensenDG_16_BSF}.
The vector case (i) can be shown analogously, and
(ii) has been proved in \cite{BramwellDGQ_12_Lhp}, see the proof of Lemma~4 around (40) there.
\end{proof}

\subsection{Trace operators for Kirchhoff--Love plates} \label{sec_KL}

Let us recall the following trace operators from \cite{FuehrerHN_19_UFK},
\begin{align*} 
   \traceGG{}:\;
   &\begin{cases}
      H^2(\cT) &\to\ \HdDiv{\cT}',\\
      \quad z    &\mapsto\ \dual{\traceGG{}(z)}{\SS}_\cS
                    := \vdual{z}{\div\Div\SS}_\cT - \vdual{\sGrad{\grad z}}{\SS}_\cT,
   \end{cases}
   \\
   \traceDD{}:\;
   &\begin{cases}
      \HdDiv{\cT} &\to\ H^2(\cT)',\\
      \qquad \SS    &\mapsto\ \dual{\traceDD{}(\SS)}{z}_\cS
                   := \dual{\traceGG{}(z)}{\SS}_\cS.
   \end{cases}
\end{align*}
Again, we also need the trace operators 
$\traceGG{\Omega}:\;H^2(\Omega)\to\HdDiv{\Omega}'$ and
$\traceDD{\Omega}:\;\HdDiv{\Omega}\to H^2(\Omega)'$.
They are defined as before, only replacing $\cT$ with $\Omega$.
To deal with the simply supported case \eqref{BCs} we introduce
\[
   \HdDivz{\Omega}
   :=
   \{\MM\in\HdDiv{\Omega};\; \dual{\traceDD{\Omega}(\MM)}{z}_\Gamma=0\
                             \forall z\in H^2(\Omega)\cap H^1_0(\Omega)\}.
\]
The following relations hold true.
\begin{prop} \label{prop_KL}
(i) If $z\in H^2(\cT)$, then
\[
   z\in H^2_0(\Omega) \quad\Leftrightarrow\quad
   \dual{\traceDD{}(\SS)}{z}_\cS = 0
   \quad\forall \SS\in \HdDiv{\Omega}
\]
and
\[
   z\in H^2(\Omega)\cap H^1_0(\Omega) \quad\Leftrightarrow\quad
   \dual{\traceDD{}(\SS)}{z}_\cS = 0
   \quad\forall \SS\in \HdDivz{\Omega}
\]
(ii) If $\SS\in \HdDiv{\cT}$, then
\[
   \SS\in \HdDiv{\Omega} \quad\Leftrightarrow\quad
   \dual{\traceGG{}(z)}{\SS}_\cS = 0
   \quad\forall z\in H^2_0(\Omega)
\]
and
\[
   \SS\in \HdDivz{\Omega} \quad\Leftrightarrow\quad
   \dual{\traceGG{}(z)}{\SS}_\cS = 0
   \quad\forall z\in H^2(\Omega)\cap H^1_0(\Omega).
\]
\end{prop}

\begin{proof}
Statements (i) are proved by \cite[Proposition~3.8]{FuehrerHN_19_UFK},
\cite[Corollary~12]{FuehrerHS_20_UFR}, and for statements (ii) we refer to
\cite[Proposition~11]{FuehrerHS_20_UFR}.
\end{proof}

\subsection{Shell-specific spaces and trace operators} \label{sec_sK}

For $T\in\cT$, we define the space
\[
   \VV(T) = \bH^1(T)\times H^2(T)\times \HDivs{T}\times \HdDiv{T} \times \LL_2^k(T)
\]
equipped with the norm (squared)
\begin{align} \label{testnorm}
   \|\vv\|_{\VV(T)}^2
   :=
   & D^{-2} \|\Cdisp\bv\|_T^2 + \|\Grad\bv - \BB z + \QQ\|_T^2
   + d^2 D^{-4} \|z\|_T^2 + d^2 \|\sGrad\grad z\|_T^2
   + \|\TT\|_T^2
   \nonumber\\
   &
   + D^2 \|\Cdispinv\Div\TT\|_T^2
   + d^{-2} \|\SS\|_T^2
   + d^{-2} D^4 \|\div\Div\SS-\BB\vdotdot\TT\|_T^2 + \cQ \|\QQ\|_T^2
\end{align}
for $\vv=(\bv,z,\TT,\SS,\QQ)$. Here, $\Cdisp$ is a positive definite, symmetric tensor of
order two (already announced in \eqref{Poincare}) and $\cQ>0$ is a constant.
Their choice will be made clear when defining the DPG scheme in Section~\ref{sec_DPG},
cf.~\eqref{ass} and \eqref{cQ}.

Analogously we define $\VV(\cT)$ as the product space with respect to $T\in\cT$,
with norm $\|\cdot\|_{\VV(\cT)}$ and inner product $\ip{\cdot}{\cdot}_{\VV(\cT)}$,
and the space $\VV(\Omega)$ with norm $\|\cdot\|_{\VV(\Omega)}$,
defined like the space $\VV(T)$ and norm $\|\cdot\|_{\VV(T)}$ just replacing $T$ with $\Omega$.
The $\VV$-spaces are used for test functions. For the solution we need the spaces
\[
   U(T) = \bH^1(T)\times H^2(T)\times \HDiv{T}\times \HdDiv{T}\quad (T\in\cT)
\]
with (squared) norm
\begin{align} \label{U_norm}
   &\|(\bu,w,\NN,\MM)\|_{U(T)}^2
   :=
     D^{-2} \|\Cdisp\bu\|_T^2 + \|\sGrad\bu + \BB w\|_T^2
   + d^2 D^{-4} \|w\|_T^2 + d^2 \|\sGrad\grad w\|_T^2
   + \|\NN\|_T^2
   \nonumber\\ &\qquad
   + \cQinv \|\skew\NN\|_T^2
   + D^2 \|\Cdispinv\Div\NN\|_T^2
   + d^{-2} \|\MM\|_T^2 + d^{-2} D^4 \|\div\Div\MM-\BB\vdotdot\NN\|_T^2
\end{align}
and the corresponding space $U(\Omega)$ with norm $\|\cdot\|_{U(\Omega)}$.
Here, $\skew\NN:=\bigl(\NN-\NN^T)/2$.
Now, to include boundary conditions, we denote by $\UBC{a}(\Omega)$ ($a\in\{s,c\}$)
the subspaces of elements $(\bu,w,\NN,\MM)\in U(\Omega)$ that satisfy \eqref{BCs}
($a=s$) or \eqref{BCc} ($a=c$).
The conditions of simple support \eqref{BCs} give the space
\begin{align*}
   \UBC{s}(\Omega)
   &=
   \bH^1_0(\Omega)\times \bigl(H^2(\Omega)\cap H^1_0(\Omega)\bigr)\times
   \HDiv{\Omega}\times \HdDivz{\Omega},
\end{align*}
cf.~\cite[Section~3.3]{FuehrerHS_20_UFR},
whereas clamped conditions \eqref{BCc} induce
\begin{align*}
   \UBC{c}(\Omega)
   &=
   \bH^1_0(\Omega)\times H^2_0(\Omega)\times \HDiv{\Omega}\times \HdDiv{\Omega}.
\end{align*}
We introduce the trace operator (``sK'' refers to \underline{s}hallow shell of \underline{K}oiter type)
\begin{align} \label{trSK}
   &\traceSK{}:\;
   \begin{cases}
      \qquad U(\Omega) &\to\qquad \VV(\cT)',\\
      (\bu,w,\NN,\MM) &\mapsto\ \dual{\traceSK{}(\bu,w,\NN,\MM)}{(\bv,z,\TT,\SS,\QQ)}_\cS
   \end{cases}
   \\
   &\text{where}\quad
   \dual{\traceSK{}(\bu,w,\NN,\MM)}{(\bv,z,\TT,\SS,\QQ)}_\cS
   \nonumber\\
   &\qquad\qquad:=
     \dual{\traceG{}(\bu)}{\TT}_\cS + \dual{\traceGG{}(w)}{\SS}_\cS
   + \dual{\traceD{}(\NN)}{\bv}_\cS - \dual{\traceDD{}(\MM)}{z}_\cS.
   \nonumber
\end{align}
It gives rise to the trace spaces
\begin{equation} \label{tr_image}
   \HSK(\cS):=\traceSK{}(U(\Omega)),\quad
   \HSKBC{a}(\cS) := \traceSK{}(\UBC{a}(\Omega))\quad (a\in\{s,c\}).
\end{equation}
A duality between the trace spaces and $\VV(\cT)$ is introduced to be consistent
with definition \eqref{trSK}. That is, for given $\tuu\in\HSK(\cS)$ and $\vv\in \VV(\cT)$,
their duality pairing is defined by
\begin{align*}
   \dual{\tuu}{\vv}_\cS := \dual{\traceSK{}(\uu_0)}{\vv}_\cS
\end{align*}
where $\uu_0\in U(\Omega)$ is such that $\traceSK{}(\uu_0) = \tuu$.

Of course, the trace spaces have four components consisting of the images of the four
trace operators $\traceG{}$, $\traceGG{}$, $\traceD{}$ and $\traceDD{}$.
But just by considering their joint action we achieve robustness of our numerical scheme.
To be specific, let us define the norms
\begin{align*}
   &\|\tuu\|_\trSK{\cS}
   :=
   \inf \Bigl\{\|\uu_0\|_{U(\Omega)};\;
               \uu_0\in U(\Omega),\ \traceSK{}(\uu_0)=\tuu\Bigr\},\\
   &\|\tuu\|_{\VV(\cT)'}
   :=
   \sup_{0\not=\vv\in \VV(\cT)}
   \frac{\dual{\tuu}{\vv}_{\cS}}{\|\vv\|_{\VV(\cT)}}
   \qquad\text{for} \quad \tuu \in \HSK(\cS).
\end{align*}

\begin{prop} \label{prop_norm_sK}
The identity
\[
   \|\tuu\|_{\VV(\cT)'} = \|\tuu\|_\trSK{\cS}
   \quad\forall \tuu\in \HSK(\cS)
\]
holds true for $d>0$.
\end{prop}

\begin{proof}
For the bound $\|\tuu\|_{\VV(\cT)'} \le \|\tuu\|_\trSK{\cS}$
we consider arbitrary elements $(\bu,w,\NN,\MM)\in U(\Omega)$ and $(\bv,z,\TT,\SS,\QQ)\in\VV(\cT)$. Then
\begin{align*}
   &\dual{\traceSK{}(\bu,w,\NN,\MM)}{(\bv,z,\TT,\SS,\QQ)}_\cS
   =
     \vdual{\bu}{\Div\TT}_\cT + \vdual{\sGrad{\bu}}{\TT}
   + \vdual{w}{\div\Div\SS}_\cT - \vdual{\sGrad\grad w}{\SS}
   \\
   &\hspace{0.4\textwidth}
   + \vdual{\NN}{\Grad\bv}_\cT + \vdual{\Div\NN}{\bv}
   + \vdual{\MM}{\sGrad\grad z}_\cT - \vdual{\div\Div\MM}{z}
   \\&\quad
   = \vdual{\bu}{\Div\TT}_\cT + \vdual{\sGrad{\bu}+\BB w}{\TT}
   + \vdual{w}{\div\Div\SS-\BB\vdotdot\TT}_\cT + \vdual{\NN}{\Grad\bv-\BB z + \QQ}_\cT
   \\&\quad
   - \vdual{\div\Div\MM-\BB\vdotdot\NN}{z} + \vdual{\Div\NN}{\bv}
   - \vdual{\sGrad\grad w}{\SS} + \vdual{\MM}{\sGrad\grad z}_\cT - \vdual{\NN}{\QQ}
   \\&\quad
   \le \|(\bu,w,\NN,\MM)\|_{U(\Omega)} \|(\bv,z,\TT,\SS,\QQ)\|_{\VV(\cT)}
\end{align*}
shows the claimed bound.
To prove the other inequality, we formally adapt the technique from
\cite{FuehrerHS_20_UFR}, see also the abstract framework in \cite[Lemma~A.10]{DemkowiczGNS_17_SDM}.
Note that, differently to the situation in \cite{FuehrerHS_20_UFR},
we are dealing with an operator where all the traces are independent and where different norms
are used in the domain of the operator and its dual. Though, as in \cite{FuehrerHS_20_UFR},
robust stability requires to combine test functions.

Given $\tuu\in\HSK(\cS)$, we first define $\vv=(\bv,z,\TT,\SS,\QQ)\in \VV(\cT)$ as the solution to
\begin{equation} \label{vv}
   \ip{\vv}{\deltavv}_{\VV(\cT)}
   =
   \dual{\tuu}{\deltavv}_\cS
   \ \forall \deltavv\in \VV(\cT)
\end{equation}
and then $\uu_0=(\bu,w,\NN,\MM)\in U(\Omega)$ as the solution to
\begin{equation} \label{uu0}
   \ip{\uu_0}{\deltauu}_{U(\Omega)}
   =
   \dual{\traceSK{}(\deltauu)}{\vv}_\cS
   \ \forall \deltauu\in U(\Omega).
\end{equation}
Here, $\ip{\cdot}{\cdot}_{U(\Omega)}$ denotes the inner product in $U(\Omega)$
(the inner product $\ip{\cdot}{\cdot}_{\VV(\cT)}$ was defined previously).
It follows that $\traceSK{}(\uu_0)=\tuu$, cf.~\cite[proof of Lemma~4]{FuehrerHS_20_UFR}.
Therefore, selecting $\deltavv=\vv$ in \eqref{vv} and $\deltauu=\uu_0$ in \eqref{uu0}, we find that
\[
   \|\vv\|_{\VV(\cT)}^2 = \dual{\tuu}{\vv}_\cS = \|\uu_0\|_{U(\Omega)}^2.
\]
This proves the inverse inequality,
\[
   \|\tuu\|_{\VV(\cT)'}
   \ge
   \frac{\dual{\tuu}{\vv}_\cS}{\|\vv\|_{\VV(\cT)}}
   =
   \|\uu_0\|_{U(\Omega)}
   \ge
   \|\tuu\|_\trSK{\cS}.
\]
\end{proof}

It remains to define a global test space with boundary conditions that will serve as the
domain of the adjoint problem. To this end we denote by
$\traceSK{\Omega}:\; U(\Omega)\to \VV(\Omega)'$ the trace operator restricted to continuous
test functions,
\begin{align*}
   &\traceSK{\Omega}:\;
   \begin{cases}
      \qquad U(\Omega) &\to\qquad \VV(\Omega)',\\
      (\bu,w,\NN,\MM) &\mapsto\ \dual{\traceSK{\Omega}(\bu,w,\NN,\MM)}{(\bv,z,\TT,\SS,\QQ)}_\Gamma
   \end{cases}
\end{align*}
with corresponding duality $\dual{\cdot}{\cdot}_\Gamma$ between
$\traceSK{\Omega}(U(\Omega))$ and $\VV(\Omega)$.
Of course, traces have four components with support on $\Gamma$,
\begin{equation} \label{trace_dec}
   \traceSK{\Omega}(\bu,w,\NN,\MM)
   = \Bigl(
      \traceG{\Omega}(\bu), \traceGG{\Omega}(w), \traceD{\Omega}(\NN), -\traceDD{\Omega}(\MM)
   \Bigr).
\end{equation}
We also need the space $\VV(\Omega)$ with imposed boundary conditions,
\[
   \VVBC{c}(\Omega) :=
   \bH^1_0(\Omega)\times H^2_0(\Omega)\times \HDivs{\Omega}\times \HdDiv{\Omega}\times \LL_2^k(\Omega)
\]
and
\[
   \VVBC{s}(\Omega) :=
   \bH^1_0(\Omega)\times \bigl(H^2(\Omega)\cap H^1_0(\Omega)\bigr)\times
   \HDivs{\Omega}\times \HdDivz{\Omega}\times \LL_2^k(\Omega).
\]

\begin{prop} \label{prop_jumps_sK}
(i) The identities
\[
   \VVBC{a}(\Omega) = \bigl\{\vv\in\VV(\Omega);\;
   \dual{\traceSK{\Omega}(\uu_0)}{\vv}_\Gamma = 0\ \forall \uu_0\in \UBC{a}(\Omega)\bigr\}
   \quad (a\in\{s,c\})
\]
hold true.

(ii) For $\vv\in \VV(\cT)$ and $a\in\{s,c\}$ we have that
\[
   \vv\in \VVBC{a}(\Omega) \quad\Leftrightarrow\quad
   \dual{\tuu}{\vv}_\cS = 0
   \quad\forall \tuu\in \HSKBC{a}(\cS).
\]
\end{prop}

\begin{proof}
Since the trace operators $\traceSK{}$ and $\traceSK{\Omega}$ have independent components as indicated in
\eqref{trace_dec}, the kernel relation
\[
   \dual{\traceSK{}(\bu,w,\NN,\MM)}{\vv}_\Gamma = 0\quad\forall (\bu,w,\NN,\MM)\in\UBC{a}(\Omega)
\]
for $\vv=(\bv,z,\TT,\SS,\QQ)\in \VVBC{a}(\cT)$ is equivalent to four separate kernel relations,
between $\bv$ and $\NN$, $z$ and $\MM$, $\TT$ and $\bu$, $\SS$ and $w$.
Statement (ii) then follows from Propositions~\ref{prop_can} and~\ref{prop_KL},
and the first statement is obtained in the same way by
simply considering the mesh $\cT=\{\Omega\}$. 
\end{proof}

\section{Variational formulation and DPG method} \label{sec_DPG}

We are now in a position to derive an ultraweak variational formulation of problem
\eqref{prob}, define our DPG scheme, and state well-posedness of the formulation
and robust quasi-optimal approximation properties of the numerical scheme.

To this end we consider problem \eqref{prob}
with general (not necessarily symmetric) tensor $\NN$, adding the symmetry condition weakly.
In this way, normal traces of $\NN$ can be approximated in the standard way.
We apply $12 d^{-2} \cC^{-1}$ and $\cC^{-1}$ to \eqref{p2} and \eqref{p3}, respectively.
Then, testing \eqref{p1} with $-z$, \eqref{p2} with $\SS$, \eqref{p3} with $\TT$, \eqref{p4} with $\bv$,
adding $\vdual{\NN}{\QQ}=0$ for $(\bv,z,\TT,\SS,\QQ)\in\VV(\cT)$,
and using the shallow Koiter shell trace operator defined before,
we obtain
\begin{align*}
    &\vdual{\bu}{\Div\TT}_\cT + \vdual{w}{\div\Div\SS-\BB\vdotdot\TT}_\cT
   + \vdual{\NN}{\cC^{-1}\TT+\Grad{\bv}-\BB z+\QQ}_\cT\\
   &+ \vdual{\MM}{12 d^{-2} \cC^{-1}\SS+\sGrad{\grad z}}_\cT
   - \dual{\traceSK{}(\bu,w,\NN,\MM)}{(\bv,z,\TT,\SS,\QQ)}_\cS
   =
   \vdual{\bp}{\bv} - \vdual{f}{z}.
\end{align*}
We define the trace as an independent variable with four algebraically independent components,
\[
   \tuu = (\tbu,\tw,\tNNn,\tMM) := \traceSK{}(\bu,w,\NN,\MM),
\]
and introduce the spaces
\[
   \UU_a :=
   \bL_2(\Omega)\times L_2(\Omega)\times \LL_2(\Omega)\times \LL_2^s(\Omega)\times \HSKBC{a}(\cS)
   \quad (a\in\{s,c\})
\]
with norm
\begin{align} \label{norm}
   \|(\bu,w,\NN,\MM,\tuu)\|_{\UU}^2
   &:=
   D^{-2} \|\Cdisp \bu\|^2 + d^2 D^{-4} \|w\|^2
   + \|\NN\|^2 + d^{-2} \|\MM\|^2 + \|\tuu\|_\trSK{\cS}^2.
\end{align}
Then our variational formulation reads:
\emph{Find $\uu:=(\bu,w,\NN,\MM,\tuu)\in \UU_a$ such that}
\begin{align} \label{VF}
   b(\uu,\vv) = L(\vv)\quad\forall\vv=(\bv,z,\TT,\SS,\QQ)\in \VV(\cT).
\end{align}
Here,
\begin{align*}
   b(\uu,\vv)
   &:=
   \vdual{\bu}{\Div\TT}_\cT + \vdual{w}{\div\Div\SS-\BB\vdotdot\TT}_\cT
   + \vdual{\NN}{\cC^{-1}\TT + \Grad{\bv}-\BB z+\QQ}_\cT\\
   &\quad + \vdual{\MM}{12 d^{-2} \cC^{-1}\SS + \sGrad{\grad z}}_\cT
   - \dual{\tuu}{\vv}_\cS
\end{align*}
and
\[
   L(\vv) := \vdual{\bp}{\bv} - \vdual{f}{z}.
\]
Of course, for the boundary conditions of simple support we select $a=s$, and
the clamped conditions require $a=c$. As mentioned earlier, our theory covers the case
of different boundary conditions (not only clamped or simply supported shells).
The needed assumption is that the symmetric, positive definite tensor $\Cdisp$ and constant $D>0$
be selected so that
\begin{equation} \label{ass}
   D^{-2}\|\Cdisp \bv\|^2 + d^2 D^{-4} \|z\|^2
   \lesssim \|\sGrad{\bv}-\BB z\|^2 +  d^2 \|\sGrad{\grad z}\|^2
   \quad\forall \vv=(\bv,z,\TT,\SS,\QQ)\in\VVBC{a}(\Omega)
\end{equation}
holds uniformly with respect to $d$, $\diam(\Omega)$ and $\BB$. Of course, in order to
control the tangential displacements $\bu$ and transverse deflection $w$ as strongly as possible,
the objective is to select $\Cdisp$ as large as possible
(e.g., in terms of its eigenvalues) and $D>0$ as small as possible, cf.~\eqref{norm}.
In the case of a clamped ($a=c$) or simply supported shell ($a=s$), the choices
\begin{equation} \label{Cdisp}
   D := \text{min\,width}(\Omega),\quad
   \Cdisp:=\mathrm{diag}(c,c)
   \quad\text{with}\quad
   c := \min\Bigl\{1,\frac {d}{D^2\|\BB\|_\infty}\Bigr\},\ 
   \min\{1,\frac 10\}:=1
\end{equation}
ensure that \eqref{ass} holds, cf.~\eqref{Poincare}.
Here, $\text{min\,width}(\Omega)$ denotes the minimum distance between
two parallel lines enclosing $\Omega$. For different boundary conditions,
giving rise to a space $U_a(\Omega)$ as the domain of the trace
operator $\traceSK{}$, cf.~\eqref{tr_image}, the space in \eqref{ass} is
\[
   \VVBC{a}(\Omega)
   :=\{\vv\in\VV(\Omega);\; \dual{\traceSK{\Omega}(\uu_0)}{\vv}_\Gamma = 0
                            \quad\forall\uu_0\in U_a(\Omega) \},
\]
cf.~Proposition~\ref{prop_jumps_sK}.
We also need Korn's inequality in the form
\begin{equation} \label{Korn}
   \|\kGrad\bv\|^2 \lesssim \|\sGrad\bv\|^2
   \quad\forall \vv=(\bv,z,\TT,\SS,\QQ)\in\VVBC{a}(\Omega)
\end{equation}
with skew-symmetric part of the gradient $\kGrad\bu:=\bigl(\Grad\bu-\Grad\bu^T\bigr)/2$
(usually written as
$\kGrad\bu=\frac 12\,\mathrm{rot}(\bu) \begin{pmatrix}0 & -1\\ 1 & 0\end{pmatrix}$
with $\mathrm{rot}(\bu):=\partial_xu_2-\partial_yu_1$ for $\bu=(u_1,u_2)^T$).
It is clear that \eqref{Korn} holds for both boundary conditions, $a=c$ and $a=s$,
cf.~\cite[(11.2.22)]{BrennerS_94_MTF}.
Furthermore, scaling reveals that the hidden constant does not depend on $\diam(\Omega)$.
The third parameter to be chosen for the norms is $\cQ$, cf.~\eqref{testnorm}, \eqref{U_norm}.
The stability analysis of the adjoint problem in Lemma~\ref{la_adj_well} below suggests to select
\begin{equation} \label{cQ}
   \cQ := \min\bigl\{1, d^2 \|\BB\|_\infty^{-2} D^{-4} \bigr\}.
\end{equation}

\begin{theorem} \label{thm_stab}
For $d\in (0,1]$ select $D$, $\Cdisp$ and $\cQ$ as in \eqref{Cdisp} and \eqref{cQ}, respectively.
For given $f\in L_2(\Omega)$ and $\bp\in\bL_2(\Omega)$ there exists
a unique solution $\uu=(\bu,w,\NN,\MM,\tuu)\in\UU_a$ to \eqref{VF} with
boundary conditions \eqref{BCs} ($a=s$) or \eqref{BCc} ($a=c$).
It is uniformly bounded,
\[
   \|\uu\|_{\UU}^2 \lesssim d^{-2} D^4 \|f\|^2 + D^2 \|\Cdispinv\bp\|^2
\]
with a hidden constant that is independent of $f$, $\bp$,
$\diam(\Omega)$, $\BB$, $\cT$, and $d\in (0,1]$.
Furthermore, $(\bu,w,\NN,\MM)\in \UBC{a}(\Omega)$ solves \eqref{prob},
and $\tuu=\traceSK{}(\bu,w,\NN,\MM)$.
\end{theorem}

A proof of this theorem is given in \S\ref{sec_pf}.
Now, to define the DPG method, we consider a (family of) discrete subspace(s)
$\UU_{a,h}\subset\UU_a$ (being set on a sequence of meshes $\cT$)
and introduce the \emph{trial-to-test operator} $\ttt:\;(\UU_s\cup\UU_c)\to \VV(\cT)$ by
\begin{align*}
   \ip{\ttt(\uu)}{\vv}_{\VV(\cT)} = b(\uu,\vv)\quad\forall\vv\in \VV(\cT).
\end{align*}
Then, the DPG method with optimal test functions for problem~\eqref{prob},
boundary conditions \eqref{BCs} ($a=s$) or \eqref{BCc} ($a=c$), and based
on the variational formulation~\eqref{VF}, is:
\emph{Find $\uu_h\in \UU_{a,h}$ such that}
\begin{align} \label{DPG}
   b(\uu_h,\ttt\deltauu) = L(\ttt\deltauu) \quad\forall\deltauu\in \UU_{a,h}.
\end{align}

The following theorem states its robust quasi-optimal best approximation.

\begin{theorem} \label{thm_DPG}
For $d\in (0,1]$ select $D$, $\Cdisp$ and $\cQ$ as in \eqref{Cdisp} and \eqref{cQ}, respectively.
Let $f\in L_2(\Omega)$, $\bp\in\bL_2(\Omega)$ and $a\in\{s,c\}$ be given.
For any finite-dimensional subspace $\UU_{a,h}\subset \UU_a$
there exists a unique solution $\uu_h\in \UU_{a,h}$ to \eqref{DPG}. It satisfies the quasi-optimal
error estimate
\[
   \|\uu-\uu_h\|_{\UU} \lesssim \|\uu-\ww\|_{\UU}
   \quad\forall\ww\in \UU_{a,h}
\]
with a hidden constant that is independent of $\cT$, $\diam(\Omega)$, $\BB$,
$d\in (0,1]$, and $\UU_{a,h}$.
\end{theorem}

A proof of this theorem is given in \S\ref{sec_pf}.

\begin{remark} \label{rem_delta}
It is obvious that uniform stability provided by Theorem~\ref{thm_stab} is not needed to
ensure quasi-optimal convergence stated in Theorem~\ref{thm_DPG}. For instance,
by the Sobolev embedding theorem, a point load $f$ at $\xx_0\in\Omega$ is well defined as a duality
$\vdual{f}{z}:=z(\xx_0)$ for $z\in H^2(\cT)$
(selecting the trace of $z$ from one side only if $\xx_0$ lies on the skeleton),
though it is not bounded uniformly with respect to $\cT$. Therefore, since
$d^2 D^{-4}\|z\|^2+ d^2\|\sGrad\grad z\|_\cT^2\le \|\vv\|_{\VV(\cT)}^2$
for any $\vv=(\bv,z,\TT,\SS,\QQ)\in\VV(\cT)$, Theorem~\ref{thm_DPG} also applies to
data $f\in H^2(\cT)'$, including point loads, see Section~\ref{pointLoad} for an example.
\end{remark}

\section{Inf-sup conditions and proofs of Theorems~\ref{thm_stab},~\ref{thm_DPG}} \label{sec_adj}

Let us formulate the adjoint problem to \eqref{prob} with
boundary conditions \eqref{BCs} ($a=s$) or \eqref{BCc} ($a=c$),
and data
$g_1\in L_2(\Omega)$, $\GG_2\in\LL_2^s(\Omega)$, $\GG_3\in\LL_2(\Omega)$, and $\bg_4\in\bL_2(\Omega)$.
\emph{Find $(\bv,z,\TT,\SS,\QQ)\in \VVBC{a}(\Omega)$ such that}
\begin{subequations} \label{adj}
\begin{alignat}{2}
    -\BB\vdotdot\TT + \div\Div\SS &= g_1  && \quad\text{in} \quad \Omega\label{ap1},\\
    12 d^{-2} \cC^{-1}\SS + \sGrad{\grad z} &= \GG_2  && \quad\text{in} \quad \Omega\label{ap2},\\
    \cC^{-1}\TT + \Grad{\bv} - \BB z + \QQ
    &= \GG_3  && \quad\text{in} \quad \Omega\label{ap3},\\
    \Div\TT &= \bg_4 && \quad\text{in} \quad \Omega\label{ap4}.
\end{alignat}
\end{subequations}
It is a fully non-homogeneous version of \eqref{prob} with transformed data, negative
in-plane displacement and Lagrangian multiplier $\QQ$.
We show that this problem is well posed.

\begin{lemma} \label{la_adj_well}
There exists a unique solution $\vv=(\bv,z,\TT,\SS,\QQ)\in \VVBC{a}(\Omega)$ to \eqref{adj}, and the bound
\[
   \|\vv\|_{\VV(\Omega)}^2 \lesssim
   d^{-2} D^4 \|g_1\|^2 + d^2 \|\GG_2\|^2 + \|\GG_3\|^2 + D^2 \|\Cdispinv \bg_4\|^2
\]
holds with a constant that is independent of $\diam(\Omega)$, $\BB$, $d\in (0,1]$ and the given data.
\end{lemma}

\begin{proof}
Decomposing $\GG_3=\GG_3^s+\GG_3^k$ into its symmetric and skew-symmetric parts
we use relations \eqref{ap2} and \eqref{ap3} to calculate $\SS$ and $\TT$,
\begin{equation} \label{S_T}
   \SS = \frac {d^2}{12}\cC(\GG_2 - \sGrad\grad z),\quad
   \TT = \cC(\GG_3^s - \sGrad\bv+\BB z).
\end{equation}
We then substitute $\SS$ and $\TT$ in \eqref{ap1}, \eqref{ap4}.
Testing \eqref{ap1} and \eqref{ap4} with $-\deltaz\in H^2(\Omega)$ and $\deltabv\in\bH^1(\Omega)$,
respectively, an application of the trace operators $\traceG{\Omega}$, $\traceGG{\Omega}$ yields
the relations
\begin{align*}
   & \vdual{\cC(\GG_3^s-\sGrad\bv+\BB z)}{\BB\deltaz}
   -
   \frac {d^2}{12} \vdual{\cC(\GG_2-\sGrad\grad z)}{\sGrad\grad\deltaz}
   -
   \dual{\traceGG{\Omega}(\deltaz)}{\SS}_\Gamma
   =
   - \vdual{g_1}{\deltaz},\\
   -&\vdual{\cC(\GG_3^s-\sGrad\bv+\BB z)}{\sGrad\deltabv}
   +
   \dual{\traceG{\Omega}(\deltabv)}{\TT}_\Gamma
   =
   \vdual{\bg_4}{\deltabv}.
\end{align*}
Adding both equations, we obtain
\begin{align} \label{coercive}
   & \vdual{\cC(\sGrad\bv-\BB z)}{\sGrad\deltabv-\BB\deltaz}
   + \frac {d^2}{12} \vdual{\cC\sGrad\grad z}{\sGrad\grad\deltaz}
   \nonumber\\
   &\qquad=
   -\vdual{g_1}{\deltaz}
   + \frac {d^2}{12} \vdual{\cC\GG_2}{\sGrad\grad\deltaz}
   + \vdual{\cC\GG_3^s}{\sGrad\deltabv-\BB\deltaz}
   + \vdual{\bg_4}{\deltabv}
\end{align}
if
\[
   \dual{\traceGG{\Omega}(\deltaz)}{\SS}_\Gamma
   =
   \dual{\traceG{\Omega}(\deltabv)}{\TT}_\Gamma
   = 0.
\]
By noting that the boundary conditions for \eqref{adj} are identical to those for
\eqref{prob} by Proposition~\ref{prop_jumps_sK}, cf.~\eqref{BCs} resp. \eqref{BCc},
it follows that assumption \eqref{ass} holds. Then, since
$\cC:\;\LL_2^s(\Omega)\to\LL_2^s(\Omega)$ is a self-adjoint, positive definite isomorphism,
we find that the bilinear form from \eqref{coercive} is coercive.
There is a unique solution
$(\bv,z)\in \bH^1_0(\Omega)\times H^2_0(\Omega)$ (clamped case, $a=c$) or
$(\bv,z)\in \bH^1_0(\Omega)\times \bigl(H^2(\Omega)\cap H^1_0(\Omega)\bigr)$ (simple support, $a=s$)
to \eqref{coercive} that satisfies
\begin{align*}
   \|\sGrad{\bv}-\BB z\|^2 + d^2 \|\sGrad{\grad z}\|^2
   \lesssim
   \ &\Bigl(d^{-2} D^4 \|g_1\|^2 + d^2 \|\GG_2\|^2
         + \|\GG_3^s\|^2 + D^2 \|\Cdispinv \bg_4\|^2\Bigr)^{1/2}\\
   &\Bigl(D^{-2} \|\Cdisp \bv\|^2 + \|\sGrad\bv-\BB z\|^2
         + d^2 D^{-4} \|z\|^2 + d^2 \|\sGrad\grad z\|^2\Bigr)^{1/2}.
\end{align*}
Now using \eqref{ass} explicitly,
we obtain
\begin{align} \label{adj_pf1}
   &D^{-2} \|\Cdisp \bv\|^2 + \|\sGrad{\bv}-\BB z\|^2
   + d^2 D^{-4} \|z\|^2 + d^2 \|\sGrad\grad z\|^2
   \nonumber\\
   &\lesssim
   d^{-2} D^4 \|g_1\|^2 + d^2 \|\GG_2\|^2 + \|\GG_3^s\|^2
   + D^2 \|\Cdispinv \bg_4\|^2.
\end{align}
Recalling that $\SS$ and $\TT$ are given by \eqref{S_T},
one sees that $(\bv,z,\TT,\SS,\QQ)\in\VVBC{a}(\Omega)$
with $\QQ:=\GG_3^k-\kGrad\bv$ solves \eqref{adj}, and is unique.
We continue to use relations \eqref{S_T} to bound
\begin{align}
   \label{adj_pf2}
   \|\TT\|^2 &\lesssim \|\GG_3^s\|^2 + \|\sGrad\bv-\BB z\|^2,\\
   \label{adj_pf3}
   d^{-2} \|\SS\|^2 &\lesssim d^2 \|\GG_2\|^2 + d^2 \|\sGrad\grad z\|^2.
\end{align}
The remaining terms are bounded by using \eqref{ap1},
\begin{equation} \label{adj_pf4}
   d^{-2} D^4 \|\div\Div\SS-\BB\vdotdot\TT\|^2 = d^{-2} D^4 \|g_1\|^2,
\end{equation}
and \eqref{ap4},
\begin{equation} \label{adj_pf5}
   D^2 \|\Cdispinv \Div\TT\|^2 = D^2 \|\Cdispinv \bg_4\|^2.
\end{equation}
Noting that
$\|\Grad\bv-\BB z + \QQ\|^2=\|\sGrad\bv-\BB z\|^2 + \|\kGrad\bv+\QQ\|^2$
and $\kGrad\bv+\QQ=\GG_3^k$, \eqref{adj_pf1} implies the bound
\begin{align} \label{adj_pf6}
   \|\Grad\bv-\BB z + \QQ\|^2 \lesssim
   d^{-2} D^4 \|g_1\|^2 + d^2 \|\GG_2\|^2 + \|\GG_3\|^2 + D^2 \|\Cdispinv \bg_4\|^2.
\end{align}
By Korn's inequality \eqref{Korn},
\begin{align*}
   \|\kGrad\bv\|^2
   \lesssim
   \|\sGrad\bv-\BB z\|^2 + \|\BB\|_\infty^2 \|z\|^2,
\end{align*}
so that relations $\QQ=\GG_3^k-\kGrad\bv$ and \eqref{adj_pf1} give
\begin{align*}
   & \|\QQ\|^2
   \le 2 \bigl(\|\kGrad\bv\|^2 + \|\GG_3^k\|^2\bigr)
   \lesssim
   \|\sGrad\bv-\BB z\|^2 + \|\BB\|_\infty^2 \|z\|^2
   + \|\GG_3^k\|^2 \nonumber\\
   &\quad\lesssim
   \max\bigl\{1,\|\BB\|_\infty^2 d^{-2} D^4 \bigr\}
   \Bigl(
      d^{-2} D^4 \|g_1\|^2 + d^2 \|\GG_2\|^2 + \|\GG_3\|^2 + D^2 \|\Cdispinv \bg_4\|^2
   \Bigr).
\end{align*}
This yields
\begin{align} \label{adj_pf7}
   \cQ \|\QQ\|^2
   \lesssim
   d^{-2} D^4 \|g_1\|^2 + d^2 \|\GG_2\|^2 + \|\GG_3\|^2 + D^2 \|\Cdispinv \bg_4\|^2,
\end{align}
cf.~\eqref{cQ} for the definition of $\cQ$.
A combination of \eqref{adj_pf1}--\eqref{adj_pf7} proves the claimed stability bound.
\end{proof}

\begin{lemma} \label{la_inj}
For $a\in\{s,c\}$ and $d>0$, the adjoint operator $B^*:\;\VV(\cT)\to \UU_a'$ is injective.
\end{lemma}

\begin{proof}
We assume that $(\bv,z,\TT,\SS,\QQ)\in \VV(\cT)$ is such that
$b(\deltabu,\deltaw,\deltaNN,\deltaMM,\deltatuu;\bv,z,\TT,\SS,\QQ)=0$ for any
$(\deltabu,\deltaw,\deltaNN,\deltaMM,\deltatuu)\in \UU_a$.
The selection of arbitrary $\deltatuu\in\HSKBC{a}(\cS)$ (and remaining functions zero)
yields $(\bv,z,\TT,\SS,\QQ)\in\VVBC{a}(\Omega)$ by Proposition~\ref{prop_jumps_sK}.
Therefore, $(\bv,z,\TT,\SS,\QQ)$ solves
\begin{align*}
   \Div\TT=0,\quad \div\Div\SS-\BB\vdotdot\TT=0,\quad
   \cC^{-1}\TT + \Grad{\bv}-\BB z + \QQ=0,\quad
   12 d^{-2}\cC^{-1}\SS + \sGrad{\grad z}=0
\end{align*}
in $\Omega$.
This is problem \eqref{adj} with homogeneous data. By Lemma~\ref{la_adj_well}, $(\bv,z,\TT,\SS,\QQ)=0$.
\end{proof}

\subsection{Proofs of Theorems~\ref{thm_stab},~\ref{thm_DPG}} \label{sec_pf}

Having all the tools at hand, the proofs are straightforward.
For the proof of Theorem~\ref{thm_stab} we check the standard conditions for
the bilinear form $b(\cdot,\cdot)$ and the functional $L$.
\begin{enumerate}
\item {\bf Boundedness of the functional.}
\begin{align*}
   &L(\vv) = \vdual{\bp}{\bv}-\vdual{f}{z}
   \\
   &\le
   \bigl(D^2 \|\Cdispinv \bp\|^2 + d^{-2} D^4 \|f\|^2\bigr)^{1/2}
   \bigl(D^{-2} \|\Cdisp \bv\|^2 + d^2 D^{-4} \|z\|^2\bigr)^{1/2}\\
   &\le
   \bigl(D^2 \|\Cdispinv \bp\|^2 + d^{-2} D^4 \|f\|^2\bigr)^{1/2}
   \|\vv\|_{\VV(\cT)}.
\end{align*}
\item {\bf Boundedness of the bilinear form.}  Taking into account Proposition~\ref{prop_norm_sK},
the bound $b(\uu,\vv)\lesssim \|\uu\|_{\UU}\|\vv\|_{\VV(\cT)}$ for any
$\uu\in\UU_a$ and $\vv\in\VV(\cT)$ is immediate.
\item {\bf Injectivity.} This is Lemma~\ref{la_inj}.
\item {\bf Inf-sup condition.} By \cite[Theorem~3.3]{CarstensenDG_16_BSF}, the inf-sup condition
\begin{align} \label{infsup}
   \sup_{0\not=(\bv,z,\TT,\SS,\QQ)\in \VV(\cT)}
   \frac {b(\bu,w,\NN,\MM,\tuu;\bv,z,\TT,\SS,\QQ)}{\|(\bv,z,\TT,\SS,\QQ)\|_{\VV(\cT)}}
   \gtrsim \|(\bu,w,\NN,\MM,\tuu)\|_{\UU}
\end{align}
for any $(\bu,w,\NN,\MM,\tuu)\in \UU_a$ ($a\in\{s,c\}$) follows from the inf-sup conditions
\begin{align}
   \label{infsup2}
   \sup_{0\not=(\bv,z,\TT,\SS,\QQ)\in \VV(\cT)}
   \frac{\dual{\tuu}{(\bv,z,\TT,\SS,\QQ)}_\cS}{\|(\bv,z,\TT,\SS,\QQ)\|_{\VV(\cT)}}
   \gtrsim
   \|\tuu\|_{\trSK{\cS}} \quad\forall\tq\in\HSKBC{a}(\cS)
\end{align}
and
\begin{align}
   \label{infsup1}
   \sup_{0\not=(\bv,z,\TT,\SS,\QQ)\in \VVBC{a}(\Omega)}
   &\frac{b(\bu,w,\NN,\MM,0;\bv,z,\TT,\SS,\QQ)}{\|(\bv,z,\TT,\SS,\QQ)\|_{\VV(\cT)}}
   \nonumber\\
   &\gtrsim \bigl(D^{-2} \|\Cdisp \bu\|^2
                + d^2 D^{-4} \|w\|^2 + \|\NN\|^2 + d^{-2}\|\MM\|^2\bigr)^{1/2}
\end{align}
for any $(\bu,w,\NN,\MM)\in \bL_2(\Omega)\times L_2(\Omega)\times \LL_2(\Omega)\times \LL_2^s(\Omega)$.

Relation \eqref{infsup2} is true (with constant $1$) by Proposition~\ref{prop_norm_sK},
and inf-sup condition \eqref{infsup1} holds due to Lemma~\ref{la_adj_well}:
Selecting data
\[
   (g_1,\GG_2,\GG_3,\bg_4):= (d^2 D^{-4} w, d^{-2} \MM, \NN, D^{-2}\Cdispt \bu)
\]
in \eqref{adj} with solution $\vv^*\in\VVBC{a}(\Omega)$, we bound
\begin{align*}
   &\sup_{0\not=\vv\in \VVBC{a}(\Omega)}
   \frac{b(\bu,w,\NN,\MM,0;\vv)}{\|\vv\|_{\VV(\cT)}}
   \ge
   \frac{\vdual{w}{g_1} + \vdual{\MM}{\GG_2} + \vdual{\NN}{\GG_3} + \vdual{\bu}{\bg_4}}
        {\|\vv^*\|_{\VV(\cT)}}
   \\
   &\gtrsim
   \frac{d^2 D^{-4} \|w\|^2 + d^{-2} \|\MM\|^2 + \|\NN\|^2 + D^{-2}\|\Cdisp\bu\|^2}
        {\bigl(d^{-2} D^4 \|g_1\|^2 + d^2\|\GG_2\|^2 + \|\GG_3\|^2 + D^2 \|\Cdispinv \bg_4\|^2\bigr)^{1/2}}
   \\
   &=
   \bigl(d^2 D^{-4} \|w\|^2 + d^{-2}\|\MM\|^2 + \|\NN\|^2 + D^{-2} \|\Cdisp \bu\|^2\bigr)^{1/2}.
\end{align*}
\end{enumerate}
It remains to note that $(\bu,w,\NN,\MM)$ solves \eqref{prob}, satisfies the selected boundary
conditions, \eqref{BCs} or \eqref{BCc}, and $\tuu=\traceSK{}(\bu,w,\NN,\MM)$.
This finishes the proof of Theorem~\ref{thm_stab}.

Now, to prove Theorem~\ref{thm_DPG}, we recall that, by design of the DPG scheme as a minimum
residual method,
\[
   \|B(\uu-\uu_h)\|_{\VV(\cT)'} = \min\{\|B(\uu-\ww)\|_{\VV(\cT)'};\; \ww\in \UU_{a,h}\}.
\]
Here, $B:\;\UU_a\to \VV(\cT)'$ is the operator induced by the bilinear form $b(\cdot,\cdot)$.
The result follows by the uniform equivalence of $\|\cdot\|_{\UU}$ and
$\|B\cdot\|_{\VV(\cT)'}$. Indeed, $\|B\cdot\|_{\VV(\cT)'}\lesssim \|\cdot\|_{\UU}$
by the uniform boundedness of $b(\cdot,\cdot)$, and the inverse estimate is
\eqref{infsup}.

\section{Numerical experiments} \label{sec_num}

In this section, we demonstrate the performance of our DPG method in different benchmark tests.
In \S\ref{ScoLo}, we consider the classical Scordelis--Lo cylindrical shell roof which
features a somewhat realistic problem setup related to the design of reinforced concrete roofs.
As the name suggests, the problem framework was adopted for benchmarking purposes of computer methods
by Scordelis and Lo in~\cite{ScordelisL_64_CAC}. Nowadays, such benchmarking is typically restricted
to the evaluation of the maximum vertical displacement in a fixed geometry, but there is some ambiguity
in the literature regarding the exact target value. This is because the solution of the problem varies
slightly depending on the specific characteristics of the assumed shell model.
Following the analytic and semi-analytic procedures presented by Fl\"{u}gge and Briassoulis in
\cite{Fluegge_60_SS,Briassoulis_02_TAB}, we have verified the analytic reference solution for
the shallow Koiter model used in our study. The problem itself features long-range angular
boundary layers emanating from the straight edges of the roof and conventional finite element methods
are subject to moderate membrane locking.

As a second test, in \S\ref{PitCourse} we consider the Pitk\"{a}ranta cylinder obstacle course
featuring a closed circular cylindrical shell which is loaded by an axially constant but angularly
varying self-balancing normal pressure.
Pure membrane and bending dominated deformations are obtained by setting the ends of the cylinder
fully clamped (\S\ref{Case1}) or leaving them free (\S\ref{Case2}), respectively,
whereas the so-called simple effect becomes activated
for simply supported sliding edges (\S\ref{Case3}), see \cite{PitkaerantaLOP_95_SDS}.
Again, the semi-analytic reference solution has been computed for each case.
Conventionally, no locking is expected for the membrane-dominated case, whereas the inextensional
bending modes present for free ends give rise to severe membrane locking.
The simple edge effect gives also rise to transverse shear locking for Naghdi type shell models.
But these are of no concern for our Koiter type model neglecting the transverse shear strains altogether.

As a final test, in \S\ref{pointLoad} we analyze the singular stress system of a shallow shell
under a normal point load, cf.~\cite{NiemiPH_07_BCP,NiemiHP_08_PLS}.
The deformation is known to vary depending on the
geometric curvature and thickness of the shell in such a way that the ``hot spot'',
which is present around the load application point for all curved geometries, spreads into
characteristic line layers along the generators of parabolic and hyperbolic shells.
By assuming doubly periodic boundary conditions, the singularity can be examined by Fourier analysis.
The related parametric effects are relatively complex but the line layers may be associated
with the usual angular layers in cylindrical shells and with the corresponding generalized edge effect
in hyperbolic geometry. For such deformations, a moderate level of membrane locking is again expected.

Before presenting our numerical results, let us comment on the definition of tensor $\Cdisp$
that influences the control of the tangential displacements, cf.~\eqref{norm}.
It must be chosen so that (continuous) test functions $\vv=(\bv,z,\TT,\SS,\QQ)\in\VVBC{a}(\Omega)$
satisfy assumption \eqref{ass}. By duality it is clear that the boundary conditions describing
$\VVBC{a}(\Omega)$ are those of the model problem to be studied.
In other words, the boundary conditions imposed on $\bv$ and $z$
are the ones for the tangential displacements $\bu$ and the vertical deflection $w$, respectively.
In all the examples of this section, the conditions on $w$ (and therefore on $z$) are such that
$\|z\|\lesssim \diam(\Omega)^2\|\sGrad\grad z\|$
by a Poincar\'e--Friedrichs estimate (cf.~\cite[(5.9.3)]{BrennerS_94_MTF}) and scaling.
Furthermore, the conditions on the tangential displacements $\bu$ (and therefore
on $\bv$) are such that rigid body motions in the plane (rigid shell modes) are eliminated.
Therefore, $\|\bv\|\lesssim \diam(\Omega) \|\sGrad\bv\|$
by Korn's inequality (cf.~\cite[(11.2.22)]{BrennerS_94_MTF}) and scaling.
Bounding $\|\sGrad\bv\|\le \|\sGrad\bv-\BB z\|+\|\BB\|_\infty \|z\|$, and selecting $D$
such that $D\simeq\diam(\Omega)$, we conclude that
\begin{align*}
   D^{-2}\|\bv\|^2
   &\lesssim \|\sGrad{\bv}-\BB z\|^2 + \|\BB\|_\infty^2 D^4 \|\sGrad{\grad z}\|^2\\
   &\le \max\{1,d^{-2} \|\BB\|_\infty^2 D^4\}
         \left(\|\sGrad{\bv}-\BB z\|^2 + d^2 \|\sGrad{\grad z}\|^2\right)
   \quad\forall \vv=(\bv,z,\TT,\SS,\QQ)\in\VVBC{a}(\Omega)
\end{align*}
for the boundary conditions specified in the examples of this section. That is, \eqref{ass} holds with
\begin{align} \label{Cdisp_num}
   D\simeq\diam(\Omega),\quad \Cdisp = \mathrm{diag}(c,c),\quad c=\min\{1,d\|\BB\|_\infty^{-1} D^{-2}\}.
\end{align}
This is a generic choice for $\Cdisp$ as mentioned in Section~\ref{sec_model}
(recall relations \eqref{Poincare} ff.) and claimed in \eqref{Cdisp} (improved for stretched domains).
In some cases below, a stronger selection is in order, and possible. We do not
study anisotropic domains here but refer to \cite{FuehrerH_21_RDM} for such an example.

\subsection{Setup}

Before presenting our numerical results, let us specify the discretization setup.
We use shape-regular triangular meshes $\cT$ on $\Omega$ with maximum mesh-size $h=\max_{T\in\cT} \diam(T)$.
In all the examples, $\Omega$ is a rectangular domain and we use an initial mesh of four triangles. 
By $\cP^k(\cT)$ we denote the space of $\cT$-piecewise polynomials of degree $\leq k\in\N_0$,
and $[\mathcal{RT}^0(\cT)]^2\subset\HDiv{\Omega}$ is the space of row-wise
lowest-order Raviart--Thomas elements. We also need the subspace
$\mathcal{rHCT}(\cT)\subset H^2(\Omega)$ consisting of the (reduced) HCT elements
(the degrees of freedom are associated with the values at vertices and
the values of gradients at vertices). For an edge $E$ let $\cP^0(E)$ be the space of constants on $E$. 
Then
\begin{align*}
  \mathcal{FHN}(\cT) := \{\MM\in\HdDiv{\Omega};\; &\nn\cdot\MM\nn|_E\in \cP^0(E),\, \nn\cdot\Div\MM + \partial_{\tangent}(\tangent\cdot\MM\nn)|_E \in \cP^0(E), \,  
  \\&\text{ for all edges } E \text{ of the triangulation }\cT\}.
\end{align*}
For details on the construction of this space (more precisely its trace space) and its approximation
properties, we refer to~\cite[Section~6]{FuehrerHN_19_UFK}.
Finally we set $\widehat\UU_{hk}$ to be the image of the trace operator $\traceSK{}$ acting on
\[
   [\cP^{k+1}(\cT)\cap H^1(\Omega)]^2\times \mathcal{HCT}(\cT)\times
   [\mathcal{RT}^0(\cT)]^2 \times \mathcal{FHN}(\cT),
\]
and use the approximation spaces
\begin{align*}
  \UU_{a,hk} :=
  \bigl([\cP^0(\cT)]^2 \times \cP^0(\cT) \times [\cP^0(\cT)]^{2\times 2} \times
   [\cP^0(\cT)]^{2\times 2}\times \widehat\UU_{hk}\bigr) \cap \UU_a, \quad k=0,1.
\end{align*}
Here, the parameter $a$ refers to the type of boundary condition. Previously we used
$a\in\{c,s\}$ for the clamped and simply supported cases. Now, this notation generically refers
to the space satisfying the required boundary conditions.
Considering an ultraweak formulation, we treat all boundary conditions like essential ones.
In the examples below, we only describe the kinematic constraints
(i.e., boundary conditions for the primal variables $w,\bu$).
The remaining conditions for $\NN$ and $\MM$, which are dual to the ones for $\bu$ and $w$,
respectively, are taken to be homogeneous.
We also stress the fact that the only difference between the cases $k=0$ and $k=1$ consists
in the increased polynomial degree for the trace of the tangential displacements.
This specific enrichment aims at alleviating the membrane locking effect.

Membrane locking is caused by the fact that the discrete approximation space for the displacements cannot
produce bending modes with vanishing membrane strains.
Depending on the exact nature of the bending mode, some or all components of the membrane strain tensor
shall vanish. The former happens, e.g., for the various boundary and interior layer modes
whereas the latter case is possible if the kinematic constraints are sufficiently weak to allow
pure bending, or isometric transformations of the middle surface.
Even though the locking mechanism for the DPG method may not be identical to that of
standard displacement based finite elements, we expect that increasing the approximation degree
of the trace variables associated with the tangential displacements shall alleviate the locking effect.
Indeed, this is what we observe in the examples below.

\begin{remark}
Static condensation of the field approximation variables is possible, and would drastically reduce the
overall degrees of freedom to the ones associated with the trace variables.
We do not consider this elimination of field variables in our experiments.

Furthermore, instead of using general (non-symmetric) approximations of $\NN$ as we do,
one can consider symmetric discretizations. This change eliminates the test variable $\QQ$ because
$(\NN,\QQ)=0$ for any $\QQ\in\LL_2^k(\cT)$ is then satisfied automatically
and since the trace operator is independent of $\QQ$, cf.~\eqref{trSK}.
\end{remark}

For the discretization of the test space we use
\begin{align*}
  \VV_h := \VV_h(\cT) :=
  [\cP^3(\cT)]^2\times \cP^3(\cT) \times [\cP^3(\cT)]^{2\times 2,\mathrm{sym}}\times
       [\cP^4(\cT)]^{2\times 2,\mathrm{sym}}\times [\cP^2(\cT)]^{2\times 2,\mathrm{skew}}
\end{align*}
where the notation ``sym'' and ``skew'' in the upper indices refers to the subsets of
symmetric and skew-symmetric tensors, respectively.
The choices of polynomial degrees for the first four components are based on the
constructions in~\cite{GopalakrishnanQ_14_APD,FuehrerH_19_FDD}
where Fortin operators are provided for a linear elasticity problem and a Kirchhoff--Love
plate bending model. They use the same discretization spaces for the field and trace variables.
Though, in our case, the construction of Fortin operators that are uniformly bounded in $d$ is
an open problem. This is due to the parameter-dependent coupling of terms in the norm
of the test space, cf.~\eqref{testnorm}.
However, our numerical examples below indicate that our choice is reasonable. 

Let
\begin{equation} \label{eta}
   \eta^2 := \sum_{T\in\cT} \eta(T)^2 :=
   \sum_{T\in\cT} \|B\uu_h-L\|_{\VV_h(T)'}^2 = \|B\uu_h-L\|_{\VV_h(\cT)'}^2
\end{equation}
denote the DPG built-in error estimator which---provided the existence of a Fortin operator---is
reliable and efficient up to an oscillation term, see~\cite{CarstensenDG_14_PEC}.
Here, $\uu_h \in \UU_{a,hk}$ is the DPG approximation and $\VV_h(T)$ denotes the discrete test space
on one element. We use this estimator to steer adaptive mesh refinements,
marking elements for refinement by the bulk criterion
\begin{align*}
  \theta \eta^2 \leq \sum_{T\in\mathcal{M}} \eta(T)^2.
\end{align*}
Here, $\mathcal{M}\subset\cT$ is a set of minimal cardinality and we choose $\theta = 1/4$.
The newest-vertex bisection rule is used to divide each marked element $T\in\mathcal{M}$
into four triangles of equal area. 
We refer to uniform refinement when each element is marked for refinement, i.e, $\mathcal{M}=\cT$.

In the following, let $\uu=(\bu,w,\NN,\MM,\widehat\uu)$ and
$\uu_h = (\bu_h,w_h,\NN_h,\MM_h,\widehat\uu_h)\in \UU_{a,hk}$ be the generic exact solution
and its DPG approximation, respectively. We denote $\bu=(u_1,u_2)$ and $\NN=(N_{ij})$, and
use the following abbreviations in the presentation of numerical results,
\begin{alignat*}{2}
  &\err(w_h) := d \|w-w_h\|,            && \err(\bu_h) := \|\Cdisp(\bu-\bu_h)\|, \\
  &\err(\MM_h) := d^{-1}\|\MM-\MM_h\|,  \qquad && \err(\NN_h) := \|\NN-\NN_h\|,
\end{alignat*}
and $\ndof$ is the dimension of $\UU_{a,hk}$. In the case of quasi-uniform meshes, $\ndof=\OO(h^{-2})$.
Note that, if $\diam(\Omega)\simeq D\simeq 1$,
then $\err(w_h)+\err(\bu_h)+\err(\MM_h)+\err(\NN_h)\lesssim \|\uu-\uu_h\|_{\UU}$,
cf.~\eqref{norm}

\subsection{Scordelis--Lo cylindrical shell roof} \label{ScoLo}

Adopting the measurement units from~\cite{MacNeal_94_FET},
the Scordelis--Lo benchmark problem can be formulated by choosing
\begin{align*}
  \Omega = (0,R)\times (0,\alpha R), \quad \BB = \begin{pmatrix} 0 & 0 \\ 0 & 1/R \end{pmatrix}
  \quad\text{with}\quad R=25,\ \alpha=2\pi/9,
\end{align*}
corresponding to a quarter of the full roof, cut at $x=0$ and $y=0$.
Recalling the employed scaling of the displacements, the forcing corresponding to a
uniformly distributed vertical load $g$ is defined as
\begin{align*}
  \bp=(p_1,p_2),\quad
  p_1 = 0, \quad p_2 = \frac g{dE} \sin\left(\frac{y}{R}\right), \quad
  f = -\frac g{dE}\cos\left(\frac{y}{R}\right).
\end{align*}
The symmetry boundary conditions at the cut lines are
\[
  u_1|_{x=0}=0, \quad \partial_{\nn} w|_{x=0} = 0, \quad u_2|_{y=0}=0, \quad \partial_{\nn} w|_{y=0}=0,
\]
and the conditions of the rigid diaphragm at the end of the roof are
\[
  u_2|_{x=R} = 0, \quad w|_{x=R} = 0.
\]
For the values
\begin{align*}
  E = 4.32 \cdot 10^8, \quad g = 90, \quad\nu = 0, \quad d=0.25,
\end{align*}
the vertical displacement at the midpoint of the free edge is evaluated as
\begin{align} \label{value}
  u_2(0,\alpha R)\sin\alpha - w(0,\alpha R)\cos\alpha \approx 0.3086.
\end{align}
For this problem, a stronger tensor $\Cdisp$ than specified in \eqref{Cdisp_num} can be taken.
Indeed, by the form of $\BB$ (vanishing curvature in $x$-direction) and the boundary condition
$v_1|_{x=0}=0$ (implied by the one for $u_1$) where $\bv=(v_1,v_2)$,
Poincar\'e--Friedrichs inequality yields
$\|v_1\|\le R\|\partial_x v_1\|\le R\|\sGrad\bv-\BB z\|$.
The $v_2$-component is controlled as in \eqref{Cdisp_num}. Therefore, \eqref{ass} holds with
$D=R$ and $\Cdisp = \mathrm{diag}(1,d/R)$, our selection for this example.

\begin{figure}
  \begin{center}
    \begin{tikzpicture}
\begin{loglogaxis}[
width=0.49\textwidth,
xlabel={$\ndof$},
grid=major,
legend entries={{\small $\eta$, $k=0$},{\small $\eta$, $k=1$}, {\small $\OO(\ndof^{-1/2})$}},
legend pos=south west,
every axis plot/.append style={ultra thick},
]
\addplot table [x=dofDPG,y=estDPG] {data/ScoLoUh0.dat};
\addplot table [x=dofDPG,y=estDPG] {data/ScoLoUh1.dat};
\addplot [black,dotted,mark=none] table [x=dofDPG,y expr={0.1*sqrt(\thisrowno{1})^(-1)}] {data/ScoLoUh0.dat};

\end{loglogaxis}
\end{tikzpicture}
\begin{tikzpicture}
\begin{loglogaxis}[
width=0.49\textwidth,
xlabel={$\ndof$},
grid=major,
legend entries={{\small rel. err., $k=0$},{\small rel. err., $k=1$}, {\small $\OO(\ndof^{-3/4})$}},
legend pos=south west,
every axis plot/.append style={ultra thick},
]
\addplot table [x=dofDPG,y=errRef] {data/ScoLoUh0.dat};
\addplot table [x=dofDPG,y=errRef] {data/ScoLoUh1.dat};
\addplot [black,dotted,mark=none] table [x=dofDPG,y expr={50*sqrt(\thisrowno{1})^(-1.5)}] {data/ScoLoUh0.dat};

\end{loglogaxis}
\end{tikzpicture}
  \end{center}
  \caption{Scordelis--Lo. Left: error estimator $\eta$,
           right: errors for reference value \eqref{value}.}\label{fig:ScoLo}
\end{figure}
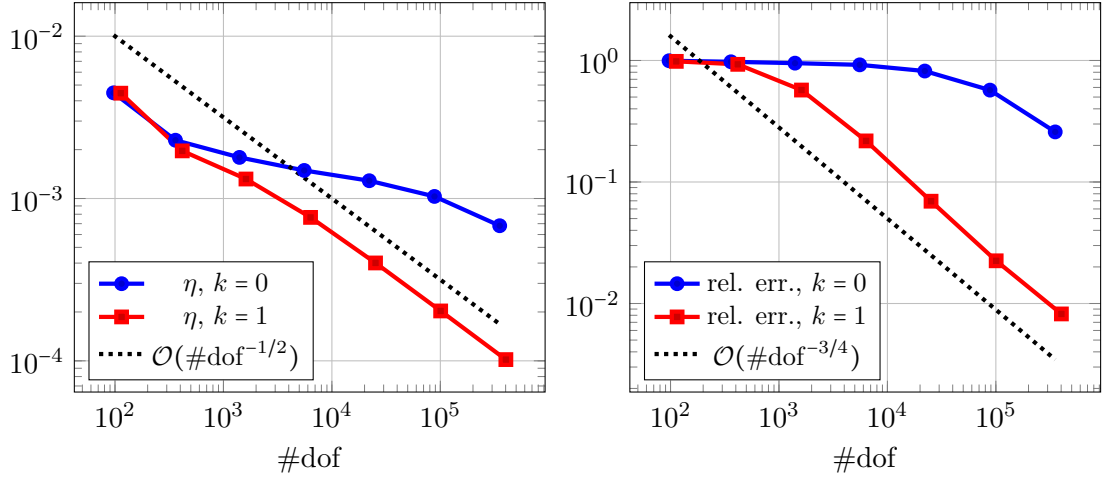

The following results are obtained with uniform mesh refinements. The left plot in
Figure~\ref{fig:ScoLo} compares the error estimators $\eta$ provided by the DPG framework, cf.~\eqref{eta}.
We observe a moderate locking effect when using the space $\UU_{a,h0}$ which is alleviated
when switching to $\UU_{a,h1}$. The right plot compares the relative errors for the
approximation of the vertical displacement \eqref{value}.
As expected we obtain better results when using $\UU_{a,h1}$.
In fact, in the case $k=1$, we observe a convergence
with rate $\OO(\ndof^{-3/4})=\OO(h^{3/2})$ after only a brief pre-asymptotic behavior. 
We note that the discrete value is computed using the piecewise constant approximations
to $w$ and $\bu$ in \eqref{value}.

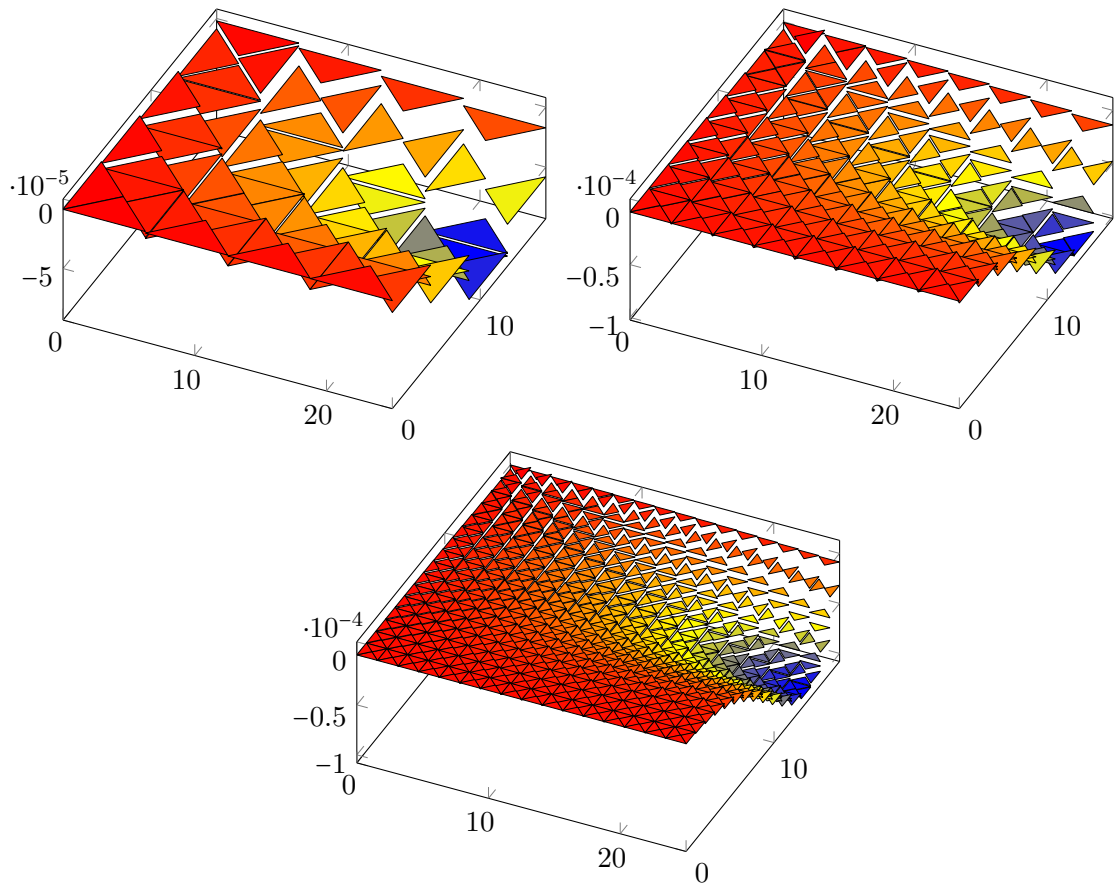
\begin{figure}
  \begin{center}
    \begin{tikzpicture}
  \begin{axis}[
width=0.5\textwidth,
view={25}{60},
]
\addplot3[patch,line width=0.2pt,faceted color=black] table{data/ScoLoN12nE64.dat};
\end{axis}
\end{tikzpicture}
\begin{tikzpicture}
  \begin{axis}[
width=0.5\textwidth,
view={25}{60},
]
\addplot3[patch,line width=0.2pt,faceted color=black] table{data/ScoLoN12nE256.dat};
\end{axis}
\end{tikzpicture}
\begin{tikzpicture}
  \begin{axis}[
width=0.5\textwidth,
view={25}{60},
]
\addplot3[patch,line width=0.2pt,faceted color=black] table{data/ScoLoN12nE1024.dat};
\end{axis}
\end{tikzpicture}
  \end{center}
  \caption{Scordelis--Lo, approximation of $N_{12}$, $k=1$. Meshes with $64$, $256$, $1024$ elements.}
  \label{fig:ScoLo:N12}
\end{figure}

An important feature of our method is its ability to approximate stress resultants robustly even on
relatively coarse meshes. This is demonstrated in Figure~\ref{fig:ScoLo:N12} showing the distribution
of the in-plane shear force approximation of $N_{12}$ on different meshes.
This quantity of interest is known to be rather difficult to approximate for some shell finite element
formulations, see, e.g., \cite{AndelfingerR_93_ETT}.

\subsection{Pitk\"{a}ranta cylindrical shell obstacle course} \label{PitCourse}

For this set of problems we choose
\begin{align*}
  \Omega = (-1,1)\times (0,\pi/4), \quad \BB = \begin{pmatrix} 0 & 0 \\ 0 & 1 \end{pmatrix},\quad
  E = 1,\quad \nu = 0,\quad f(x,y) = \cos(2y)
\end{align*}
and throughout the boundary conditions
\begin{align} \label{BC_cyl}
  w|_{y=\pi/4} = 0, \quad 
  \partial_{\nn}w|_{y=0} = 0, \quad 
  u_1|_{y=\pi/4} =0, \quad
  u_2|_{y=0}=0.
\end{align}
They represent the symmetry of the problem on the full cylinder,
restricted in our experiments to angles between $0$ and $\pi/4$.
Additional kinematic constraints will be specified for the individual cases,
considered in the following subsections. Throughout we select the parameter $D=1$.

\subsubsection{Membrane state: cylinder with clamped ends} \label{Case1}

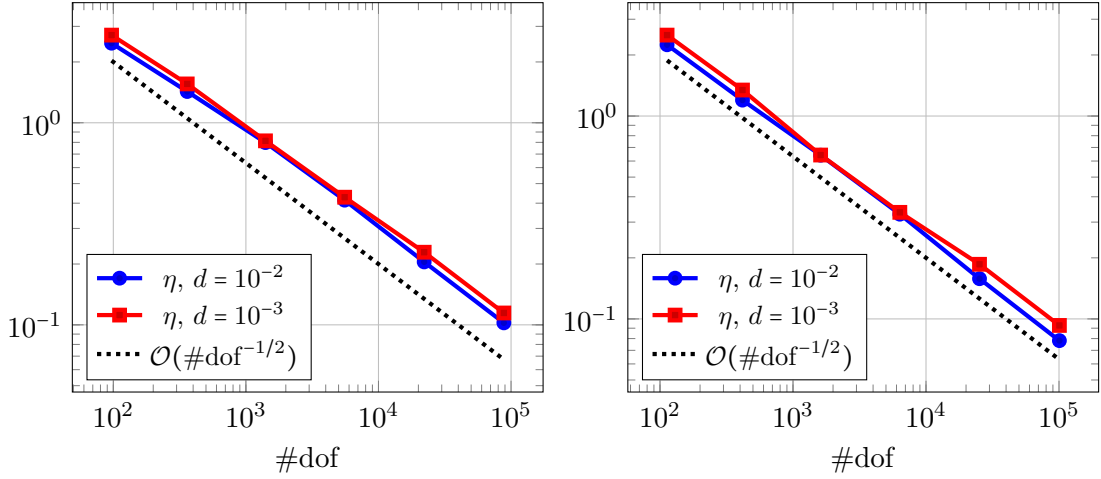
\begin{figure}
  \begin{center}
    \begin{tikzpicture}
\begin{loglogaxis}[
width=0.49\textwidth,
xlabel={$\ndof$},
grid=major,
legend entries={{\small $\eta$, $d=10^{-2}$},{\small $\eta$, $d=10^{-3}$}, {\small $\OO(\ndof^{-1/2})$}},
legend pos=south west,
every axis plot/.append style={ultra thick},
]
\addplot table [x=dofDPG,y=estDPG] {data/Case1Uh0_d2.dat};
\addplot table [x=dofDPG,y=estDPG] {data/Case1Uh0_d3.dat};
\addplot [black,dotted,mark=none] table [x=dofDPG,y expr={20*sqrt(\thisrowno{1})^(-1)}] {data/Case1Uh0_d2.dat};

\end{loglogaxis}
\end{tikzpicture}
\begin{tikzpicture}
\begin{loglogaxis}[
width=0.49\textwidth,
xlabel={$\ndof$},
grid=major,
legend entries={{\small $\eta$, $d=10^{-2}$},{\small $\eta$, $d=10^{-3}$}, {\small $\OO(\ndof^{-1/2})$}},
legend pos=south west,
every axis plot/.append style={ultra thick},
]
\addplot table [x=dofDPG,y=estDPG] {data/Case1Uh1_d2.dat};
\addplot table [x=dofDPG,y=estDPG] {data/Case1Uh1_d3.dat};
\addplot [black,dotted,mark=none] table [x=dofDPG,y expr={20*sqrt(\thisrowno{1})^(-1)}] {data/Case1Uh1_d2.dat};

\end{loglogaxis}
\end{tikzpicture}
  \end{center}
  \caption{Clamped cylinder, error estimator $\eta$.
           Left: $k=0$, right: $k=1$.}\label{fig:Case1}
\end{figure}

Additionally to \eqref{BC_cyl}, all displacements and the rotation are set to vanish along the curved ends,
\begin{align*}
  w|_{x=\pm 1} = 0,\quad \partial_\nn w|_{x=\pm 1}=0, \quad \bu|_{x=\pm 1} = 0.
\end{align*}
By the arguments given previously in \S\ref{ScoLo}, and noting that $\diam(\Omega)\simeq 1$,
\eqref{ass} is satisfied with our selection $\Cdisp = \mathrm{diag}(1,d)$.
This case is membrane dominated and no relevant locking effect is expected.
This is reflected by the results in Figure~\ref{fig:Case1} which show
the built-in error estimator for a sequence of uniform meshes.
The values $d=10^{-2},10^{-3}$ are considered, on the left for $k=0$ and on the right for $k=1$.

\subsubsection{Inextensional state: cylinder with free ends} \label{Case2}

\begin{figure}
  \begin{center}
    \begin{tikzpicture}
\begin{loglogaxis}[
width=0.49\textwidth,
cycle list name=color list,
cycle multiindex* list={
mark list*\nextlist
color list\nextlist
},
xlabel={$\ndof$},
grid=major,
legend entries={{\small $\eta$},{\small $\err(w_h)$}, {\small $\err(\bu_h)$}, {\small $\err(\MM_h)$}, {\small
$\err(\NN_h)$}, {\small $\OO(\ndof^{-1/2})$}},
legend style={at={(0.5,1.05)},anchor=south},
legend columns=3, 
        legend style={
            /tikz/column 2/.style={
                column sep=5pt,
            }},
every axis plot/.append style={ultra thick},
]
\addplot table [x=dofDPG,y=estDPG] {data/Case2Uh0_d2.dat};
\addplot table [x=dofDPG,y=errW] {data/Case2Uh0_d2.dat};
\addplot table [x=dofDPG,y=errU] {data/Case2Uh0_d2.dat};
\addplot table [x=dofDPG,y=errM] {data/Case2Uh0_d2.dat};
\addplot table [x=dofDPG,y=errN] {data/Case2Uh0_d2.dat};
\addplot [black,dotted,mark=none] table [x=dofDPG,y expr={100*sqrt(\thisrowno{1})^(-1)}] {data/Case2Uh0_d2.dat};

\end{loglogaxis}
\end{tikzpicture}
\begin{tikzpicture}
\begin{loglogaxis}[
width=0.49\textwidth,
cycle list name=color list,
cycle multiindex* list={
mark list*\nextlist
color list\nextlist
},
xlabel={$\ndof$},
grid=major,
legend entries={{\small $\eta$},{\small $\err(w_h)$}, {\small $\err(\bu_h)$}, {\small $\err(\MM_h)$}, {\small
$\err(\NN_h)$}, {\small $\OO(\ndof^{-1/2})$}},
legend style={at={(0.5,1.05)},anchor=south},
legend columns=3, 
        legend style={
            /tikz/column 2/.style={
                column sep=5pt,
            }},
every axis plot/.append style={ultra thick},
]
\addplot table [x=dofDPG,y=estDPG] {data/Case2Uh1_d2.dat};
\addplot table [x=dofDPG,y=errW] {data/Case2Uh1_d2.dat};
\addplot table [x=dofDPG,y=errU] {data/Case2Uh1_d2.dat};
\addplot table [x=dofDPG,y=errM] {data/Case2Uh1_d2.dat};
\addplot table [x=dofDPG,y=errN] {data/Case2Uh1_d2.dat};
\addplot [black,dotted,mark=none] table [x=dofDPG,y expr={100*sqrt(\thisrowno{1})^(-1)}] {data/Case2Uh1_d2.dat};

\end{loglogaxis}
\end{tikzpicture}
  \end{center}
  \caption{Free cylinder, errors and estimator $\eta$, $d=10^{-2}$.
           Left: $k=0$, right: $k=1$.}\label{fig:Case2}
\end{figure}
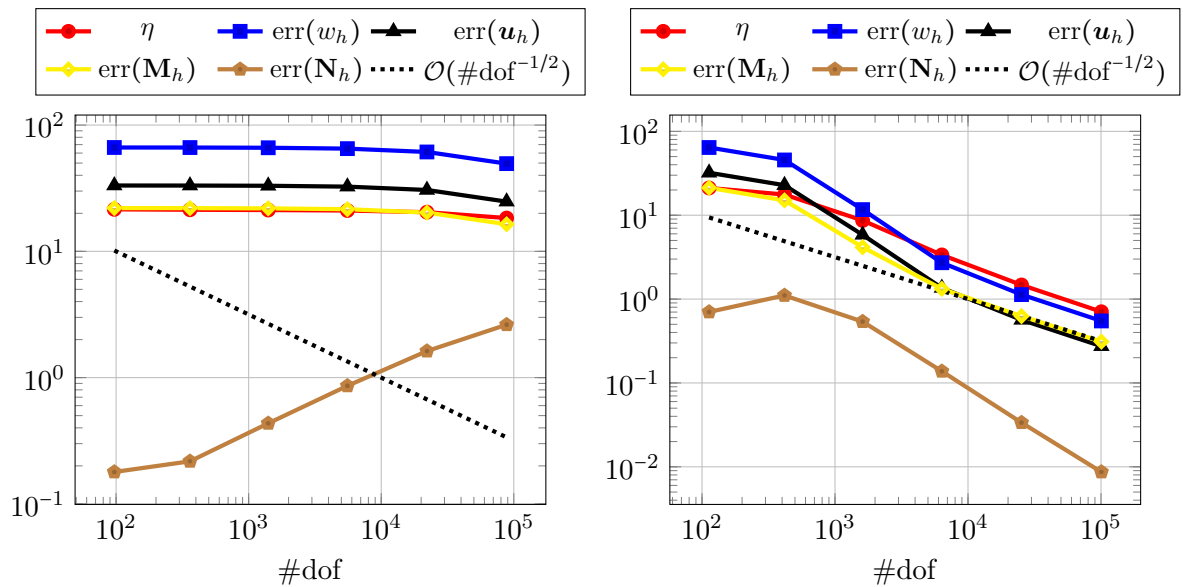

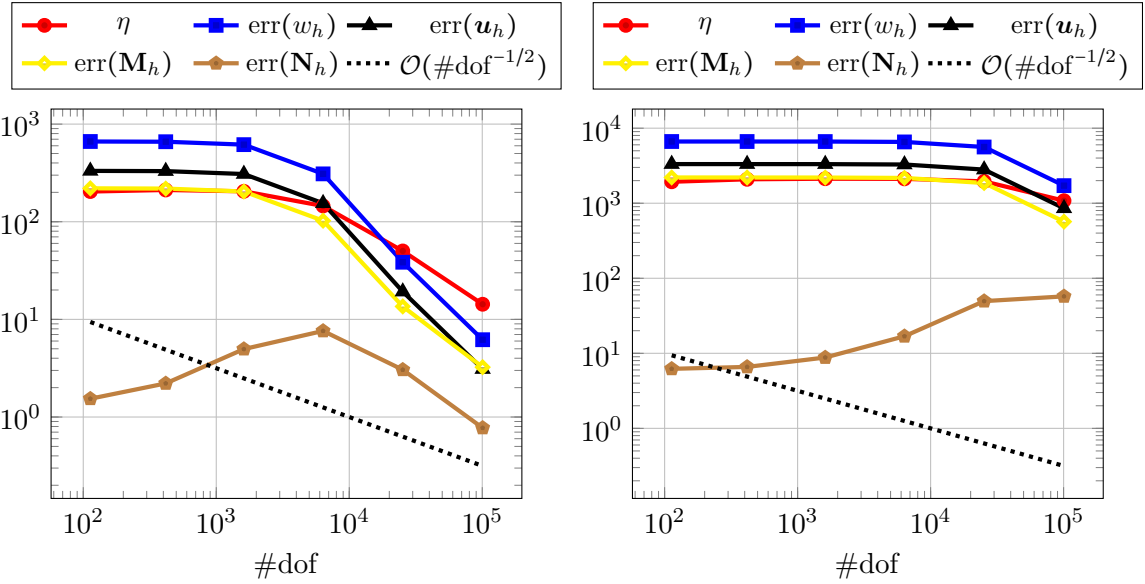
\begin{figure}
  \begin{center}
    \begin{tikzpicture}
\begin{loglogaxis}[
width=0.49\textwidth,
cycle list name=color list,
cycle multiindex* list={
mark list*\nextlist
color list\nextlist
},
xlabel={$\ndof$},
grid=major,
legend entries={{\small $\eta$},{\small $\err(w_h)$}, {\small $\err(\bu_h)$}, {\small $\err(\MM_h)$}, {\small
$\err(\NN_h)$}, {\small $\OO(\ndof^{-1/2})$}},
legend style={at={(0.5,1.05)},anchor=south},
legend columns=3, 
        legend style={
            /tikz/column 2/.style={
                column sep=5pt,
            }},
every axis plot/.append style={ultra thick},
]
\addplot table [x=dofDPG,y=estDPG] {data/Case2Uh1_d3.dat};
\addplot table [x=dofDPG,y=errW] {data/Case2Uh1_d3.dat};
\addplot table [x=dofDPG,y=errU] {data/Case2Uh1_d3.dat};
\addplot table [x=dofDPG,y=errM] {data/Case2Uh1_d3.dat};
\addplot table [x=dofDPG,y=errN] {data/Case2Uh1_d3.dat};
\addplot [black,dotted,mark=none] table [x=dofDPG,y expr={100*sqrt(\thisrowno{1})^(-1)}] {data/Case2Uh1_d3.dat};

\end{loglogaxis}
\end{tikzpicture}
\begin{tikzpicture}
\begin{loglogaxis}[
width=0.49\textwidth,
cycle list name=color list,
cycle multiindex* list={
mark list*\nextlist
color list\nextlist
},
xlabel={$\ndof$},
grid=major,
legend entries={{\small $\eta$},{\small $\err(w_h)$}, {\small $\err(\bu_h)$}, {\small $\err(\MM_h)$}, {\small
$\err(\NN_h)$}, {\small $\OO(\ndof^{-1/2})$}},
legend style={at={(0.5,1.05)},anchor=south},
legend columns=3, 
        legend style={
            /tikz/column 2/.style={
                column sep=5pt,
            }},
every axis plot/.append style={ultra thick},
]
\addplot table [x=dofDPG,y=estDPG] {data/Case2Uh1_d4.dat};
\addplot table [x=dofDPG,y=errW] {data/Case2Uh1_d4.dat};
\addplot table [x=dofDPG,y=errU] {data/Case2Uh1_d4.dat};
\addplot table [x=dofDPG,y=errM] {data/Case2Uh1_d4.dat};
\addplot table [x=dofDPG,y=errN] {data/Case2Uh1_d4.dat};
\addplot [black,dotted,mark=none] table [x=dofDPG,y expr={100*sqrt(\thisrowno{1})^(-1)}] {data/Case2Uh1_d4.dat};

\end{loglogaxis}
\end{tikzpicture}
  \end{center}
  \caption{Free cylinder, errors and estimator $\eta$, $k=1$.
           Left: $d=10^{-3}$, right: $d=10^{-4}$.}\label{fig:Case2smalld}
\end{figure}

In this case we only impose the symmetry boundary conditions \eqref{BC_cyl} and leave
the curved ends free of kinematic constraints. This gives rise to the inextensional mode
\begin{align*}
  w(x,y) = \frac 3{4d^2} \cos 2y, \quad \bu = -\frac 3{8d^2} \begin{pmatrix} 0 \\ \sin 2y\end{pmatrix}
\end{align*}
with vanishing membrane strains. Due to the lack of additional constraints
we select $\Cdisp = \mathrm{diag}(d,d)$, as in \eqref{Cdisp_num}.
It is known that this example causes severe locking when using some standard finite element methods,
cf.~\cite{PitkaerantaLOP_95_SDS}. We consider uniform meshes.
Figure~\ref{fig:Case2} shows the errors and estimator for the moderate value $d=10^{-2}$.
We observe severe locking when using the space $\UU_{a,hk}$ with $k=0$ (left plot).
It is almost eliminated when increasing the polynomial degree for the approximation
of the trace of the tangential displacements by one, $k=1$ (right plot).
As expected, locking re-appears when considering smaller values of $d$.
Figure~\ref{fig:Case2smalld} shows the results for $d=10^{-3}$ (left side) and $d=10^{-4}$ (right side),
in both cases using the enriched trace space, $k=1$. For $d=10^{-3}$ acceptable rates are observed
after a pre-asymptotic phase whereas $d=10^{-4}$ causes severe locking.

\subsubsection{Intermediate state: cylinder with simple sliding support at the ends} \label{Case3}

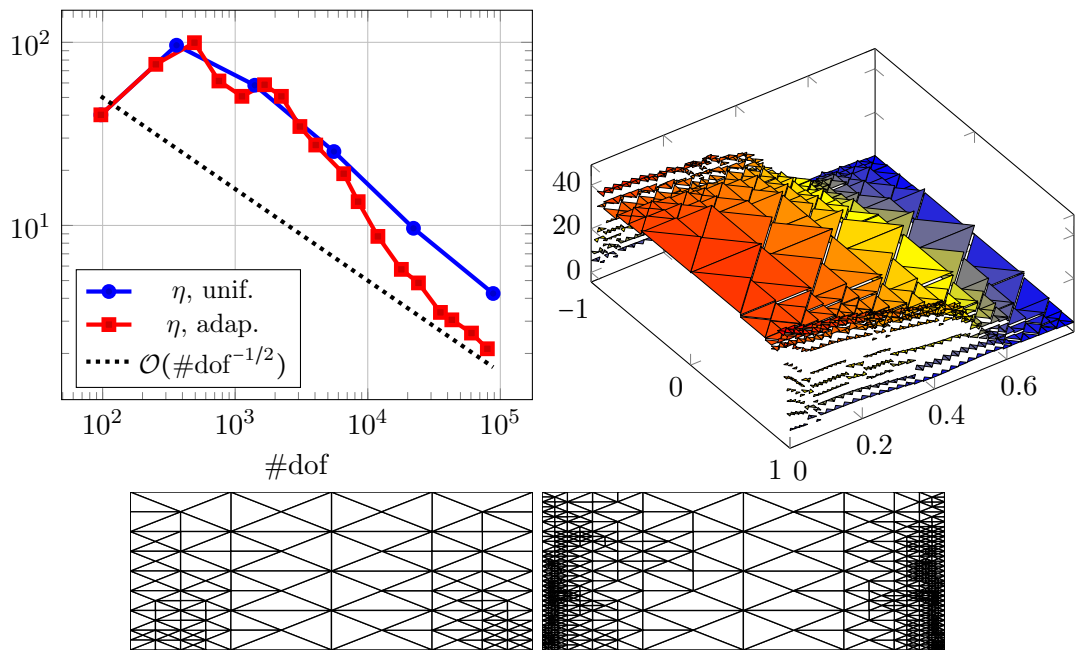
\begin{figure}
  \begin{center}
    \begin{tikzpicture}
\begin{loglogaxis}[
width=0.49\textwidth,
xlabel={$\ndof$},
grid=major,
legend entries={{\small $\eta$, unif.},{\small $\eta$, adap.}, {\small $\OO(\ndof^{-1/2})$}},
legend pos=south west,
every axis plot/.append style={ultra thick},
]
\addplot table [x=dofDPG,y=estDPG] {data/Case3Uh0_d3_unif.dat};
\addplot table [x=dofDPG,y=estDPG] {data/Case3Uh0_d3_adap.dat};
\addplot [black,dotted,mark=none] table [x=dofDPG,y expr={500*sqrt(\thisrowno{1})^(-1)}] {data/Case3Uh0_d3_unif.dat};

\end{loglogaxis}
\end{tikzpicture}
\begin{tikzpicture}
  \begin{axis}[
width=0.5\textwidth,
view={55}{60},
]
\addplot3[patch,line width=0.2pt,faceted color=black] table{data/Case3SolW.dat};
\end{axis}
\end{tikzpicture}

\begin{tikzpicture}
\begin{axis}[hide axis,
width=0.5\textwidth,
    axis equal,
]

\addplot[patch,color=white,
faceted color = black, line width = 0.5pt,
patch table ={data/Case3ele_nE186.dat}] file{data/Case3coo_nE186.dat};
\end{axis}
\end{tikzpicture}
\begin{tikzpicture}
\begin{axis}[hide axis,
width=0.5\textwidth,
    axis equal,
]

\addplot[patch,color=white,
faceted color = black, line width = 0.5pt,
patch table ={data/Case3ele.dat}] file{data/Case3coo.dat};
\end{axis}
\end{tikzpicture}
  \end{center}
  \caption{Cylinder with simple sliding support, $d=10^{-3}$.
           Top left: error estimator $\eta$ for uniform and adaptive refinements.
           Top right: transverse deflection $w_h$ on adaptively refined mesh, $\#\cT=1115$.
           Bottom: adaptively refined meshes, $\#\cT=186$ (left), $\#\cT=1115$ (right).
           Equal axis scaling in bottom row.}\label{fig:Case3}
\end{figure}

We impose boundary conditions \eqref{BC_cyl} as before,
and consider cylinder ends with simple sliding support,
\begin{align*}
  w|_{x=\pm1} = 0.
\end{align*}
As in the case of free ends, the choice $\Cdisp=\mathrm{diag}(d,d)$ implies that \eqref{ass} holds.
In this case the analytic expression of the solution is cumbersome. It can be computed reliably
with computer algebra software that utilizes arbitrary-precision arithmetic.
However, it is well known that for this problem boundary layers occur at $x=\pm 1$,
cf.~\cite{PitkaerantaLOP_95_SDS}, and the effect of locking does not play an important role here. 
We consider $d=10^{-3}$ and compare the results for uniform and adaptive mesh refinements.
In Figure~\ref{fig:Case3} we present our results when using $\UU_{a,h0}$
(the results for $\UU_{a,h1}$ are basically the same and are not shown). They show that
the adaptive variant nicely identifies the layer zones and reaches sooner the range of optimal convergence.

\subsection{Shallow shell subject to a concentrated load} \label{pointLoad}

We now consider shallow shells under a normal point load (delta distribution) at $(0,0)$,
i.e., $f=\delta_{(0,0)}$ and $\bp=0$, cf.~Remark~\ref{rem_delta}.
We study cases with elliptic, parabolic and hyperbolic middle surface:
\begin{align*}
  \Omega = (-1,1)\times (-1,1),\quad
  \BB_{\mathrm{ell}} = \begin{pmatrix} 1 & 0 \\ 0 & 1 \end{pmatrix}, \quad
  \BB_{\mathrm{par}} = \begin{pmatrix} 0 & 0 \\ 0 & 1 \end{pmatrix}, \quad
  \BB_{\mathrm{hyp}} = \begin{pmatrix} 0 & 1 \\ 1 & 0 \end{pmatrix},
\end{align*}
select $E=1$, $\nu=0$, $D=1$, and consider $d=10^{-2}$ for all the examples in this subsection.
For the actual computations we replace the delta distribution by a discrete version
$\delta_{h,(0,0)}$ which is supported on only one element $T_0\in\cT$ such that $(0,0)$
is a vertex of $T_0$ and $(\delta_{h,(0,0)},z) = z|_{T_0}(0,0)$ for all $z\in\cP^3(\cT)$.

For the elliptic and parabolic cases we impose the boundary conditions
\begin{align*}
  w|_{\partial \Omega} = 0, \quad u_1|_{y=\pm 1} = 0, \quad u_2|_{x=\pm 1}=0
\end{align*}
corresponding to the expansions
\begin{equation} \label{expansions}
  \begin{aligned}
  w(x,y) &= \sum_{m,n=1}^\infty W_{mn} \cos (m-1/2)\pi x \,\cos(n-1/2)\pi y, \\
  u_1(x,y) &= \sum_{m,n=1}^\infty \alpha_{mn} \sin(m-1/2)\pi x \,\cos(n-1/2)\pi y, \\
  u_2(x,y) &= \sum_{m,n=1}^\infty \beta_{mn} \cos(m-1/2)\pi x \,\sin(n-1/2)\pi y,
\end{aligned}
\end{equation}
whereas for the hyperbolic case we impose
\begin{align*}
  w|_{\partial \Omega} = 0, \quad u_1|_{x=\pm 1} = 0, \quad u_2|_{y=\pm 1}=0
\end{align*}
corresponding to
\begin{equation} \label{rotatedexpansions}
  \begin{aligned}
  w(x,y) &= \sum_{m,n=1}^\infty W_{mn} \cos (m-1/2)\pi x \,\cos(n-1/2)\pi y, \\
  u_1(x,y) &= \sum_{m,n=1}^\infty \alpha_{mn} \cos(m-1/2)\pi x \,\sin(n-1/2)\pi y, \\
  u_2(x,y) &= \sum_{m,n=1}^\infty \beta_{mn} \sin(m-1/2)\pi x \,\cos(n-1/2)\pi y.
\end{aligned}
\end{equation}

\subsubsection{Elliptic case}\label{pointLoadElliptic}

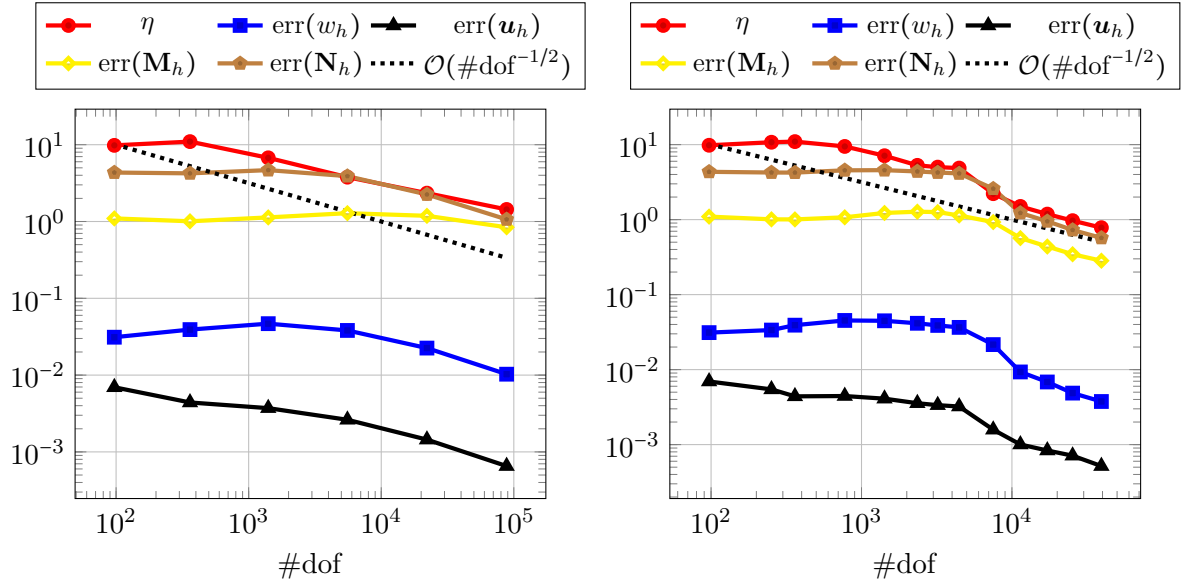
\begin{figure}
  \begin{center}
    \begin{tikzpicture}
\begin{loglogaxis}[
width=0.49\textwidth,
cycle list name=color list,
cycle multiindex* list={
mark list*\nextlist
color list\nextlist
},
xlabel={$\ndof$},
grid=major,
legend entries={{\small $\eta$},{\small $\err(w_h)$}, {\small $\err(\bu_h)$}, {\small $\err(\MM_h)$}, {\small
$\err(\NN_h)$}, {\small $\OO(\ndof^{-1/2})$}},
legend style={at={(0.5,1.05)},anchor=south},
legend columns=3, 
        legend style={
            /tikz/column 2/.style={
                column sep=5pt,
            }},
every axis plot/.append style={ultra thick},
]
\addplot table [x=dofDPG,y=estDPG] {data/PointLoadEllipticUh0_d2_unif.dat};
\addplot table [x=dofDPG,y=errW] {data/PointLoadEllipticUh0_d2_unif.dat};
\addplot table [x=dofDPG,y=errU] {data/PointLoadEllipticUh0_d2_unif.dat};
\addplot table [x=dofDPG,y=errM] {data/PointLoadEllipticUh0_d2_unif.dat};
\addplot table [x=dofDPG,y=errN] {data/PointLoadEllipticUh0_d2_unif.dat};
\addplot [black,dotted,mark=none] table [x=dofDPG,y expr={100*sqrt(\thisrowno{1})^(-1)}] {data/PointLoadEllipticUh0_d2_unif.dat};

\end{loglogaxis}
\end{tikzpicture}
\begin{tikzpicture}
\begin{loglogaxis}[
width=0.49\textwidth,
cycle list name=color list,
cycle multiindex* list={
mark list*\nextlist
color list\nextlist
},
xlabel={$\ndof$},
grid=major,
legend entries={{\small $\eta$},{\small $\err(w_h)$}, {\small $\err(\bu_h)$}, {\small $\err(\MM_h)$}, {\small
$\err(\NN_h)$}, {\small $\OO(\ndof^{-1/2})$}},
legend style={at={(0.5,1.05)},anchor=south},
legend columns=3, 
        legend style={
            /tikz/column 2/.style={
                column sep=5pt,
            }},
every axis plot/.append style={ultra thick},
]
\addplot table [x=dofDPG,y=estDPG] {data/PointLoadEllipticUh0_d2_adap.dat};
\addplot table [x=dofDPG,y=errW] {data/PointLoadEllipticUh0_d2_adap.dat};
\addplot table [x=dofDPG,y=errU] {data/PointLoadEllipticUh0_d2_adap.dat};
\addplot table [x=dofDPG,y=errM] {data/PointLoadEllipticUh0_d2_adap.dat};
\addplot table [x=dofDPG,y=errN] {data/PointLoadEllipticUh0_d2_adap.dat};
\addplot [black,dotted,mark=none] table [x=dofDPG,y expr={100*sqrt(\thisrowno{1})^(-1)}] {data/PointLoadEllipticUh0_d2_adap.dat};

\end{loglogaxis}
\end{tikzpicture}
  \end{center}
  \caption{Elliptic shell with point load, $d=10^{-2}$, $k=0$, errors and estimator $\eta$.
           Left: uniform meshes, right: adaptively refined meshes.}\label{fig:pointLoadElliptic}
\end{figure}

\begin{figure}
  \begin{center}
    \begin{tikzpicture}
  \begin{axis}[
width=0.5\textwidth,
legend entries={{\small exact},{\small unif.}, {\small adap.}},
xlabel={$x$},
]
\addplot[ultra thick,black,dotted] table [x=x,y=N11] {data/PointLoadEllipticN11exactUh0_d2_unif.dat};
\addplot[thick,red] table [x=x,y=N11] {data/PointLoadEllipticN11Uh0_d2_unif.dat};
\addplot[thick,blue] table [x=x,y=N11] {data/PointLoadEllipticN11Uh0_d2_adap.dat};
\end{axis}
\end{tikzpicture}
  \end{center}
  \caption{Elliptic shell with point load, $d=10^{-2}$, $k=0$.
           Exact solution $N_{11}$ along $y=0$ and its approximations
           with uniform mesh ($4096$ triangles) and adaptively refined mesh ($1828$ triangles).}
  \label{fig:pointLoadEllipticN11}
\end{figure}
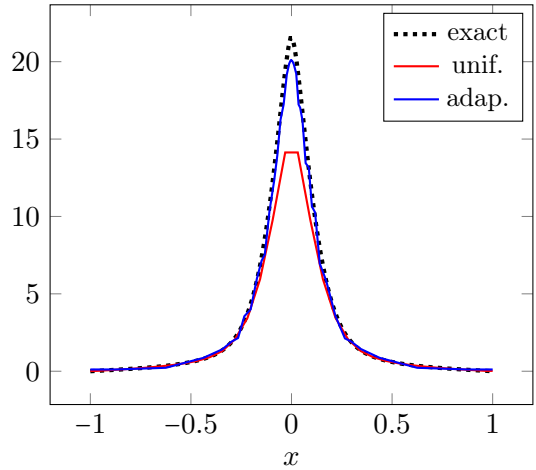

In the elliptic case, the Fourier coefficients associated with \eqref{expansions} are
\begin{align*}
  W_{mn} &= \frac{12}{d^2(M^2+N^2)^2+12},\quad
  \alpha_{mn}=\frac{-12M}{d^2(M^2+N^2)^3 + 12(M^2+N^2)}, \\
  \beta_{mn}&= \frac{-12N}{d^2(M^2+N^2)^3 + 12(M^2+N^2)} 
  \qquad\text{where}\quad M: = (m-1/2)\pi,\ N:=(n-1/2)\pi.
\end{align*}
According to \eqref{Cdisp_num} we choose $\Cdisp=\mathrm{diag}(d,d)$.
In Figure~\ref{fig:pointLoadElliptic} we compare the errors for the field variables and the
estimators for uniformly and adaptively refined meshes, in both cases using $k=0$.
For $k=1$ we get comparable results (not shown). 
We observe a rather long pre-asymptotic phase for the errors $\err(\MM_h)$ and $\err(N_h)$
which is reduced by adaptivity.
Figure~\ref{fig:pointLoadEllipticN11} shows the solution $N_{11}$ along the line $y=0$
together with its approximations ($k=0$) on the finest meshes of uniform and adaptive refinements.
We observe that the maximum value of the approximation on the adaptive mesh provides a good
approximation of the exact value.

\subsubsection{Parabolic case}\label{pointLoadParabolic}

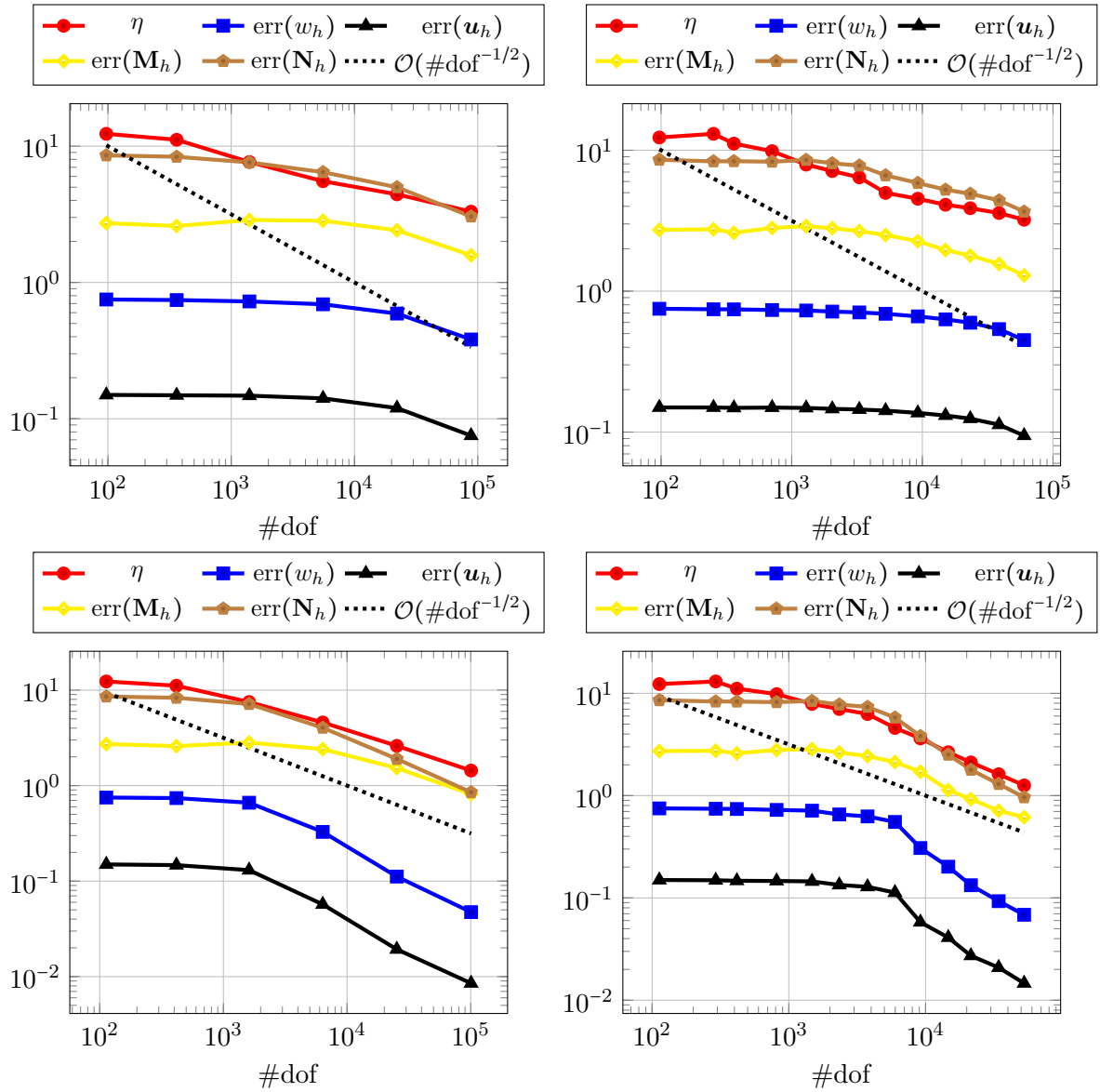
\begin{figure}
  \begin{center}
    \begin{tikzpicture}
\begin{loglogaxis}[
width=0.49\textwidth,
cycle list name=color list,
cycle multiindex* list={
mark list*\nextlist
color list\nextlist
},
xlabel={$\ndof$},
grid=major,
legend entries={{\small $\eta$},{\small $\err(w_h)$}, {\small $\err(\bu_h)$}, {\small $\err(\MM_h)$}, {\small $\err(\NN_h)$}, {\small $\OO(\ndof^{-1/2})$}},
legend style={at={(0.5,1.05)},anchor=south},
legend columns=3, 
        legend style={
            /tikz/column 2/.style={
                column sep=5pt,
            }},
every axis plot/.append style={ultra thick},
]
\addplot table [x=dofDPG,y=estDPG] {data/PointLoadParabolicUh0_d2_unif.dat};
\addplot table [x=dofDPG,y=errW] {data/PointLoadParabolicUh0_d2_unif.dat};
\addplot table [x=dofDPG,y=errU] {data/PointLoadParabolicUh0_d2_unif.dat};
\addplot table [x=dofDPG,y=errM] {data/PointLoadParabolicUh0_d2_unif.dat};
\addplot table [x=dofDPG,y=errN] {data/PointLoadParabolicUh0_d2_unif.dat};
\addplot [black,dotted,mark=none] table [x=dofDPG,y expr={100*sqrt(\thisrowno{1})^(-1)}] {data/PointLoadParabolicUh0_d2_unif.dat};

\end{loglogaxis}
\end{tikzpicture}
\begin{tikzpicture}
\begin{loglogaxis}[
width=0.49\textwidth,
cycle list name=color list,
cycle multiindex* list={
mark list*\nextlist
color list\nextlist
},
xlabel={$\ndof$},
grid=major,
legend entries={{\small $\eta$},{\small $\err(w_h)$}, {\small $\err(\bu_h)$}, {\small $\err(\MM_h)$}, {\small $\err(\NN_h)$}, {\small $\OO(\ndof^{-1/2})$}},
legend style={at={(0.5,1.05)},anchor=south},
legend columns=3, 
        legend style={
            /tikz/column 2/.style={
                column sep=5pt,
            }},
every axis plot/.append style={ultra thick},
]
\addplot table [x=dofDPG,y=estDPG] {data/PointLoadParabolicUh0_d2_adap.dat};
\addplot table [x=dofDPG,y=errW] {data/PointLoadParabolicUh0_d2_adap.dat};
\addplot table [x=dofDPG,y=errU] {data/PointLoadParabolicUh0_d2_adap.dat};
\addplot table [x=dofDPG,y=errM] {data/PointLoadParabolicUh0_d2_adap.dat};
\addplot table [x=dofDPG,y=errN] {data/PointLoadParabolicUh0_d2_adap.dat};
\addplot [black,dotted,mark=none] table [x=dofDPG,y expr={100*sqrt(\thisrowno{1})^(-1)}] {data/PointLoadParabolicUh0_d2_adap.dat};

\end{loglogaxis}
\end{tikzpicture}
\begin{tikzpicture}
\begin{loglogaxis}[
width=0.49\textwidth,
cycle list name=color list,
cycle multiindex* list={
mark list*\nextlist
color list\nextlist
},
xlabel={$\ndof$},
grid=major,
legend entries={{\small $\eta$},{\small $\err(w_h)$}, {\small $\err(\bu_h)$}, {\small $\err(\MM_h)$}, {\small $\err(\NN_h)$}, {\small $\OO(\ndof^{-1/2})$}},
legend style={at={(0.5,1.05)},anchor=south},
legend columns=3, 
        legend style={
            /tikz/column 2/.style={
                column sep=5pt,
            }},
every axis plot/.append style={ultra thick},
]
\addplot table [x=dofDPG,y=estDPG] {data/PointLoadParabolicUh1_d2_unif.dat};
\addplot table [x=dofDPG,y=errW] {data/PointLoadParabolicUh1_d2_unif.dat};
\addplot table [x=dofDPG,y=errU] {data/PointLoadParabolicUh1_d2_unif.dat};
\addplot table [x=dofDPG,y=errM] {data/PointLoadParabolicUh1_d2_unif.dat};
\addplot table [x=dofDPG,y=errN] {data/PointLoadParabolicUh1_d2_unif.dat};
\addplot [black,dotted,mark=none] table [x=dofDPG,y expr={100*sqrt(\thisrowno{1})^(-1)}] {data/PointLoadParabolicUh1_d2_unif.dat};

\end{loglogaxis}
\end{tikzpicture}
\begin{tikzpicture}
\begin{loglogaxis}[
width=0.49\textwidth,
cycle list name=color list,
cycle multiindex* list={
mark list*\nextlist
color list\nextlist
},
xlabel={$\ndof$},
grid=major,
legend entries={{\small $\eta$},{\small $\err(w_h)$}, {\small $\err(\bu_h)$}, {\small $\err(\MM_h)$}, {\small $\err(\NN_h)$}, {\small $\OO(\ndof^{-1/2})$}},
legend style={at={(0.5,1.05)},anchor=south},
legend columns=3, 
        legend style={
            /tikz/column 2/.style={
                column sep=5pt,
            }},
every axis plot/.append style={ultra thick},
]
\addplot table [x=dofDPG,y=estDPG] {data/PointLoadParabolicUh1_d2_adap.dat};
\addplot table [x=dofDPG,y=errW] {data/PointLoadParabolicUh1_d2_adap.dat};
\addplot table [x=dofDPG,y=errU] {data/PointLoadParabolicUh1_d2_adap.dat};
\addplot table [x=dofDPG,y=errM] {data/PointLoadParabolicUh1_d2_adap.dat};
\addplot table [x=dofDPG,y=errN] {data/PointLoadParabolicUh1_d2_adap.dat};
\addplot [black,dotted,mark=none] table [x=dofDPG,y expr={100*sqrt(\thisrowno{1})^(-1)}] {data/PointLoadParabolicUh1_d2_adap.dat};

\end{loglogaxis}
\end{tikzpicture}
  \end{center}
  \caption{Parabolic shell with point load, $d=10^{-2}$, errors and estimator $\eta$.
           Left: uniform meshes, right: adaptively refined meshes,
           top: $k=0$, bottom: $k=1$.}\label{fig:pointLoadParabolic}
\end{figure}

\begin{figure}
  \begin{center}
    \begin{tikzpicture}
  \begin{axis}[
width=0.5\textwidth,
legend entries={{\small exact},{\small unif.}, {\small adap.}},
xlabel={$y$},
]
\addplot[ultra thick,black,dotted] table [x=y,y=N22] {data/PointLoadParabolicN22exactUh1_d2_unif.dat};
\addplot[thick,red] table [x=y,y=N22] {data/PointLoadParabolicN22Uh1_d2_unif.dat};
\addplot[thick,blue] table [x=y,y=N22] {data/PointLoadParabolicN22Uh1_d2_adap.dat};
\end{axis}
\end{tikzpicture}
  \end{center}
  \caption{Parabolic shell with point load, $d=10^{-2}$, $k=1$.
           Exact solution $N_{22}$ along $x=0$ and its approximations
           with uniform mesh ($4096$ triangles) and adaptively refined mesh ($2139$ triangles).}
  \label{fig:pointLoadParabolicN22}
\end{figure}
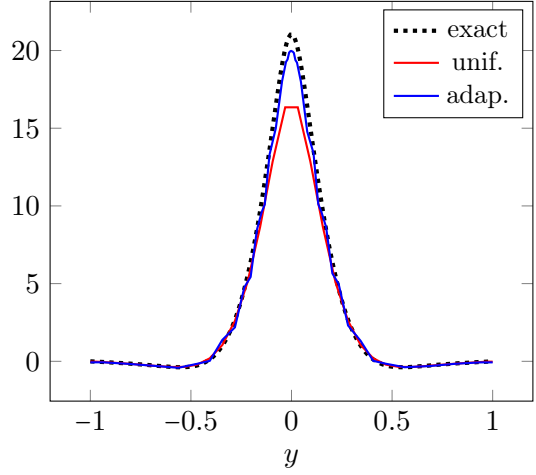

In the parabolic case, the Fourier coefficients associated with \eqref{expansions} are
\begin{align*}
  W_{mn} &= \frac{12(M^2+N^2)^2}{d^2(M^2+N^2)^4+12M^4},\quad
  \alpha_{mn} = \frac{12MN^2}{d^2(M^2+N^2)^4+12M^4},\\
  \beta_{mn} &= \frac{-24M^2N + 12N^3}{d^2(M^2+N^2)^4 + 12M^4}
  \qquad\text{where}\quad
  M: = (m-1/2)\pi,\ N:=(n-1/2)\pi.
\end{align*}
As in the elliptic case we set $\Cdisp=\mathrm{diag}(d,d)$.
In Figure~\ref{fig:pointLoadParabolic} we compare the errors of the field variables
and the estimator for uniform and adaptively refined meshes, in both cases for $k=0$ and $k=1$. 
Contrary to the elliptic case, we observe stronger locking when using $k=0$.
In Figure~\ref{fig:pointLoadParabolicN22} we compare the solution $N_{22}$ along $x=0$
with its approximations on the finest uniform and adaptively refined meshes, using $k=1$.
Again, the adaptive version provides a good approximation of the maximum value of $N_{22}$.

\subsubsection{Hyperbolic case}\label{pointLoadHyperbolic}

\begin{figure}
  \begin{center}
    \begin{tikzpicture}
\begin{loglogaxis}[
width=0.49\textwidth,
cycle list name=color list,
cycle multiindex* list={
mark list*\nextlist
color list\nextlist
},
xlabel={$\ndof$},
grid=major,
legend entries={{\small $\eta$},{\small $\err(w_h)$}, {\small $\err(\bu_h)$}, {\small $\err(\MM_h)$}, {\small
$\err(\NN_h)$}, {\small $\OO(\ndof^{-1/2})$}},
legend style={at={(0.5,1.05)},anchor=south},
legend columns=3, 
        legend style={
            /tikz/column 2/.style={
                column sep=5pt,
            }},
every axis plot/.append style={ultra thick},
]
\addplot table [x=dofDPG,y=estDPG] {data/PointLoadHyperbolicUh1_d2_unif.dat};
\addplot table [x=dofDPG,y=errW] {data/PointLoadHyperbolicUh1_d2_unif.dat};
\addplot table [x=dofDPG,y=errU] {data/PointLoadHyperbolicUh1_d2_unif.dat};
\addplot table [x=dofDPG,y=errM] {data/PointLoadHyperbolicUh1_d2_unif.dat};
\addplot table [x=dofDPG,y=errN] {data/PointLoadHyperbolicUh1_d2_unif.dat};
\addplot [black,dotted,mark=none] table [x=dofDPG,y expr={100*sqrt(\thisrowno{1})^(-1)}] {data/PointLoadHyperbolicUh1_d2_unif.dat};

\end{loglogaxis}
\end{tikzpicture}
\begin{tikzpicture}
\begin{loglogaxis}[
width=0.49\textwidth,
cycle list name=color list,
cycle multiindex* list={
mark list*\nextlist
color list\nextlist
},
xlabel={$\ndof$},
grid=major,
legend entries={{\small $\eta$},{\small $\err(w_h)$}, {\small $\err(\bu_h)$}, {\small $\err(\MM_h)$}, {\small
$\err(\NN_h)$}, {\small $\OO(\ndof^{-1/2})$}},
legend style={at={(0.5,1.05)},anchor=south},
legend columns=3, 
        legend style={
            /tikz/column 2/.style={
                column sep=5pt,
            }},
every axis plot/.append style={ultra thick},
]
\addplot table [x=dofDPG,y=estDPG] {data/PointLoadHyperbolicUh1_d2_adap.dat};
\addplot table [x=dofDPG,y=errW] {data/PointLoadHyperbolicUh1_d2_adap.dat};
\addplot table [x=dofDPG,y=errU] {data/PointLoadHyperbolicUh1_d2_adap.dat};
\addplot table [x=dofDPG,y=errM] {data/PointLoadHyperbolicUh1_d2_adap.dat};
\addplot table [x=dofDPG,y=errN] {data/PointLoadHyperbolicUh1_d2_adap.dat};
\addplot [black,dotted,mark=none] table [x=dofDPG,y expr={100*sqrt(\thisrowno{1})^(-1)}] {data/PointLoadHyperbolicUh1_d2_adap.dat};

\end{loglogaxis}
\end{tikzpicture}
  \end{center}
  \caption{Hyperbolic shell with point load, $d=10^{-2}$, $k=1$, errors and estimator $\eta$.
           Left: uniform meshes, right: adaptively refined meshes.}\label{fig:pointLoadHyperbolic}
\end{figure}
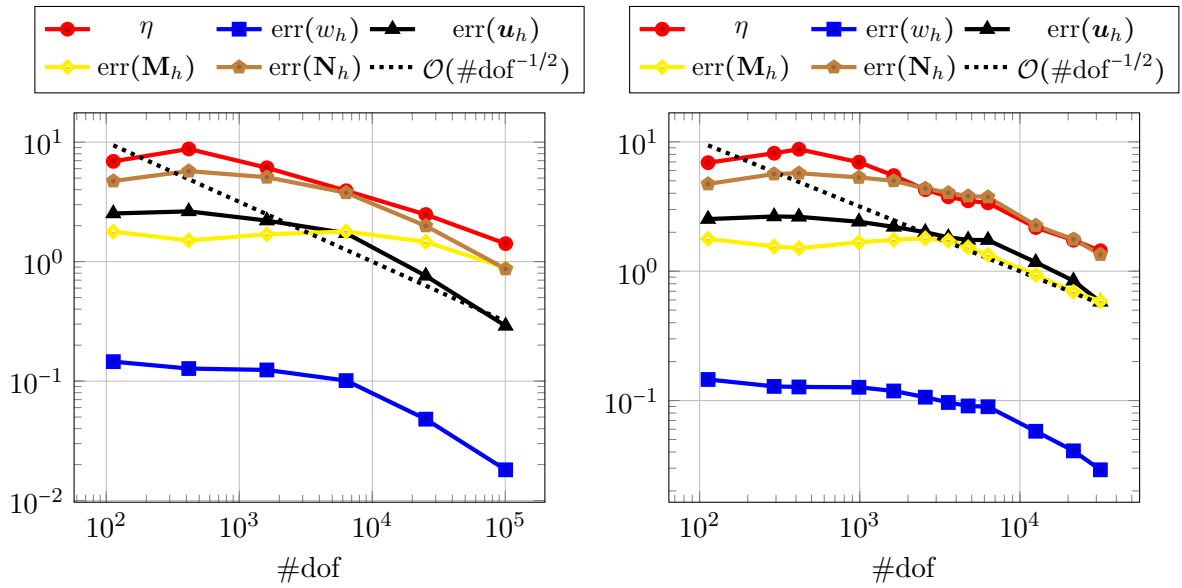

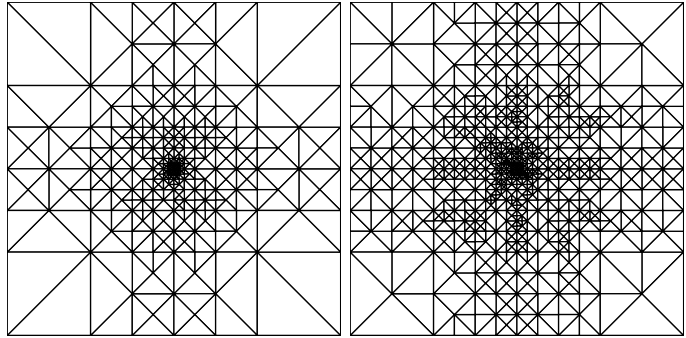
\begin{figure}
  \begin{center}
    \begin{tikzpicture}
\begin{axis}[hide axis,
width=0.5\textwidth,
    axis equal,
]

\addplot[patch,color=white,
faceted color = black, line width = 0.5pt,
patch table ={data/HyperbolicMesh1ele.dat}] file{data/HyperbolicMesh1coo.dat};
\end{axis}
\end{tikzpicture}
\begin{tikzpicture}
\begin{axis}[hide axis,
width=0.5\textwidth,
    axis equal,
]

\addplot[patch,color=white,
faceted color = black, line width = 0.5pt,
patch table ={data/HyperbolicMesh2ele.dat}] file{data/HyperbolicMesh2coo.dat};
\end{axis}
\end{tikzpicture}
  \end{center}
  \caption{Hyperbolic shell with point load, $d=10^{-2}$, $k=1$. Adaptively generated meshes,
           left: $510$ triangles, right: $1294$ triangles.}
  \label{fig:pointLoadHyperbolic:meshes}
\end{figure}

\begin{figure}
  \begin{center}
\begin{tikzpicture}
  \begin{groupplot}[view={0}{90},
    point meta min=-15,
    point meta max=75,
    group style = {group size = 2 by 1
    }]
    \nextgroupplot[colorbar horizontal,
    every colorbar/.append style={width=
        2*\pgfkeysvalueof{/pgfplots/parent axis width}+
        \pgfkeysvalueof{/pgfplots/group/horizontal sep}}]
        \addplot3[patch,hide axis] table{data/HyperbolicSol.dat};
    \nextgroupplot
    \addplot3[patch,shader=interp] table{data/HyperbolicSolExact.dat};
  \end{groupplot}
\end{tikzpicture}
  \end{center}
  \caption{Hyperbolic shell with point load, $d=10^{-2}$, transverse deflection.
           Left: approximation $w_h$ on adaptively refined mesh ($k=1$), right: exact solution $w$.}
  \label{fig:pointLoadSol}
\end{figure}
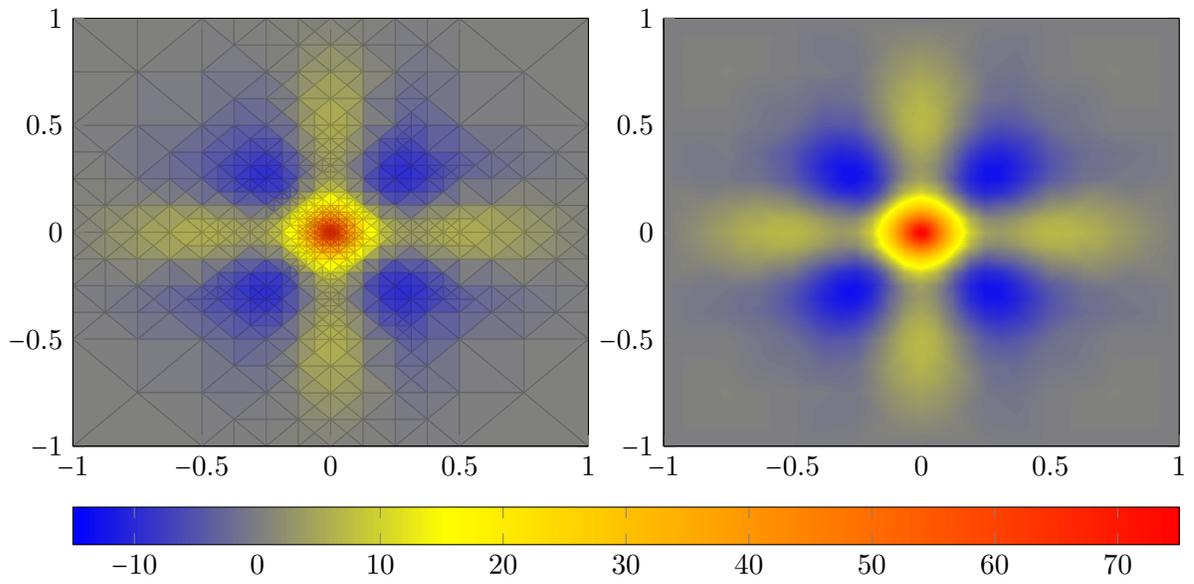

Finally, we consider the hyperbolic shell with point load. In this case,
the Fourier coefficients associated with \eqref{rotatedexpansions} are
\begin{align*}
  W_{mn} &= \frac{12(M^2+N^2)^2}{d^2(M^2+N^2)^4+48M^2N^2},\quad
  \alpha_{mn} = \frac{-24N^3}{d^2(M^2+N^2)^4+48M^2N^2},\\
  \beta_{mn} &= \frac{-24M^3}{d^2(M^2+N^2)^4+48M^2N^2}\qquad\text{where}\quad
  M: = (m-1/2)\pi,\ N:=(n-1/2)\pi.
\end{align*}
By the structure of $\BB$ (vanishing diagonal) and the boundary conditions
$v_1|_{x=\pm 1}=v_2|_{y=\pm1}=0$ (corresponding to the conditions for $\bu$),
we conclude that $\|\bv\|\le 2 \|\sGrad\bv-\BB z\|$ so that \eqref{ass} holds with
$\Cdisp = \mathrm{diag}(1,1)$.
We select $d=10^{-2}$, use the space $\UU_{a,hk}$ with $k=1$, and consider
uniform and adaptive mesh refinements, as before with parameter $\theta=1/4$.
The DPG estimator $\eta$ and errors are presented in Figure~\ref{fig:pointLoadHyperbolic}. 
We observe that, although the adaptive version does not necessarily lead to much smaller errors
in all the variables, optimal rates for the field variables are seen for fewer degrees of freedom.
Figure~\ref{fig:pointLoadHyperbolic:meshes} shows two adaptively refined meshes.
They show that the adaptive algorithm is able to localize distinctive features of the exact solution.
In Figure~\ref{fig:pointLoadSol} we compare an approximation $w_h$ of the transverse deflection
$w$ generated by the adaptive algorithm (on the left) with the exact solution $w$ (on the right).
Here, $w$ is calculated by cutting off its Fourier series after the first 10.000 terms
($m,n=1,\ldots,100$).

\FloatBarrier

\bibliographystyle{siam}
\bibliography{bib}
\end{document}